\newcount\Comments  
\documentclass[letterpaper,10pt]{article}

\usepackage[latin1]{inputenc}
\usepackage{amsmath,amstext,amsbsy,amsopn,amscd,amsxtra,upref,amssymb}
\usepackage{amsfonts,bm}
\usepackage{amsthm}
\usepackage{graphicx,cite}
\usepackage{algorithmicx}
\usepackage[]{mcode}
\usepackage{xspace}
\usepackage{algorithm}
\usepackage{algpseudocode}
\usepackage{sidecap}
\usepackage{caption}
\usepackage{lineno,hyperref}
\usepackage{cleveref}
\usepackage{fullpage}

\DeclareCaptionFormat{mylst}{\hrule#1#2#3}
\crefname{lstlisting}{listing}{listings}
\Crefname{lstlisting}{Listing}{Listings}

\usepackage{matlab-prettifier}
\usepackage{todonotes}
\usepackage{soul}
\reversemarginpar

\definecolor{ao}{rgb}{0.0, 0.5, 0.0}
\definecolor{mygreen}{rgb}{0.01, 0.75, 0.24}

\newtheorem{thm}{Theorem}[section]
\newtheorem{remark}[thm]{Remark}

\newtheorem{lem}[thm]{Lemma}
\newtheorem{cor}[thm]{Corollary}

\newenvironment{keywords}%
   {\begin{trivlist}\item[]{\bfseries\sffamily Keywords:}\ }
   {\end{trivlist}}

\ifpdf
\hypersetup{
  pdftitle={Stability of discrete empirical interpolation and gappy proper orthogonal decomposition with randomized and deterministic sampling points}, 
  pdfauthor={Benjamin Peherstorfer and Zlatko Drmac and Serkan Gugercin} 
}
\fi

\newcommand{\bfd}{\boldsymbol d}
\newcommand{\bfx}{\boldsymbol x}
\newcommand{\bfX}{\boldsymbol X}

\newcommand{\bfomega}{\boldsymbol \omega}
\newcommand{\bfxi}{\boldsymbol \xi}
\newcommand{\bff}{\boldsymbol f}
\newcommand{\bfg}{\boldsymbol g}

\newcommand{\Dcal}{\mathcal{D}}

\newcommand{\Ucal}{\mathcal{U}}

\newcommand{\bfsigma}{\boldsymbol \sigma}

\newcommand{\bfc}{\boldsymbol c}
\newcommand{\bfe}{\boldsymbol e}

\newcommand{\bfr}{\boldsymbol r}
\newcommand{\bfu}{\boldsymbol u}
\newcommand{\bfv}{\boldsymbol v}
\newcommand{\bfw}{\boldsymbol w}
\newcommand{\bfz}{\boldsymbol z}

\newcommand{\bfA}{\boldsymbol A}

\newcommand{\bfC}{\boldsymbol C}
\newcommand{\bfD}{\boldsymbol D}

\newcommand{\bfI}{\boldsymbol I}

\newcommand{\bfM}{\boldsymbol M}

\newcommand{\bfP}{\boldsymbol P}

\newcommand{\bfU}{\boldsymbol U}
\newcommand{\bfV}{\boldsymbol V}

\newcommand{\bfW}{\boldsymbol W}

\newcommand{\stilde}[1]  {\tilde{#1\hspace{0.75ex}}\hspace{-0.5ex}}
\newcommand{\shat}[1]  {\hat{#1\hspace{0.75ex}}\hspace{-0.5ex}}
\newcommand{\rowu}{\boldsymbol{u}}

\newcommand{\bfSigma}{\boldsymbol \Sigma}
\newcommand{\bfLambda}{\boldsymbol\Lambda}
\newcommand{\bfPhi}{\boldsymbol \Phi}
\newcommand{\bfOmega}{\boldsymbol \Omega}
\newcommand{\bfPsi}{\boldsymbol \Psi}

\newcommand{\bfepsilon}{\boldsymbol{\epsilon}}

\newcommand{\nh}{N} 
\newcommand{\nr}{r} 
\newcommand{\nd}{n} 
\newcommand{\np}{m} 

\newcommand{\pdim}{d} 

\newcommand{\kibitz}[2]{\ifnum\Comments=1\textcolor{#1}{#2}\fi}

\newcommand{\rpod}{\textsf{POD}\xspace}
\newcommand{\gpod}{\textsf{GappyPOD}\xspace}
\newcommand{\deim}{\textsf{DEIM}\xspace}
\newcommand{\ddeim}{(\textsf{D})\textsf{EIM}\xspace}
\newcommand{\qdeim}{\textsf{QDEIM}\xspace}

\newcommand{\odeime}{\textsf{GappyPOD+E}\xspace}

\newcommand{\odeimd}{\textsf{GappyPOD+D}\xspace}
\newcommand{\odeimr}{\textsf{GappyPOD+R}\xspace}
\newcommand{\odeiml}{\textsf{GappyPOD+L}\xspace}

\newcommand{\Id}{\boldsymbol I}

\newcommand{\bla}{\boldsymbol{\langle}}
\newcommand{\bra}{\boldsymbol{\rangle}}

\def\cool#1{#1}

		\author{
		Benjamin Peherstorfer\thanks{Courant Institute of Mathematical Sciences, New York University, New York, NY 10012. The work of Peherstorfer is supported in parts by the Air Force Center of Excellence on Multi-Fidelity Modeling of Rocket Combustor Dynamics, Award Number FA9550-17-1-0195.} \and
		Zlatko Drma\v{c}\thanks{Faculty of Science, Department of Mathematics, University of Zagreb,
		Bijeni\v{c}ka 30, 10000 Zagreb, Croatia. The work of Drma\v{c} is supported in parts by the Croatian Science Foundation through Grant IP-2019-04-6268 (Randomized low rank algorithms and applications to parameter dependent problems).} \and Serkan Gugercin\thanks{Department of Mathematics and Computational Modeling and Data Analytics Division, Academy of Integrated Science, Virginia Tech,
		Blacksburg, VA 24061-0123. The work of Gugercin is supported in parts by NSF
through Grants  DMS-1522616 and
DMS-1819110.}}
\title{Stability of discrete empirical interpolation and gappy proper orthogonal decomposition with randomized and deterministic sampling points}
\begin{document}
	\maketitle
	\begin{abstract}
		This work investigates the stability of (discrete) empirical interpolation for nonlinear model reduction and state field approximation from measurements.
		Empirical interpolation derives approximations from a few samples (measurements) via interpolation in low-dimensional spaces. It has been observed that empirical interpolation can become unstable if the samples are perturbed due to, e.g., noise, turbulence, and numerical inaccuracies. The main contribution of this work is a probabilistic analysis that shows that stable approximations are obtained if samples are randomized and if more samples than dimensions of the low-dimensional spaces are used. Oversampling, i.e., taking more sampling points than dimensions of the low-dimensional spaces, leads to approximations via regression and is known under the name of gappy proper orthogonal decomposition. Building on the insights of the probabilistic analysis, a deterministic sampling strategy is presented that aims to achieve lower approximation errors with fewer points than randomized sampling by taking information about the low-dimensional spaces into account. Numerical results of reconstructing velocity fields from noisy measurements of combustion processes and model reduction in the presence of noise demonstrate the instability of empirical interpolation and the stability of gappy proper orthogonal decomposition with oversampling.
	\end{abstract}




\section{Introduction}
\label{sec:Introduction}
Model reduction seeks to construct reduced systems that provide accurate approximations of the solutions of large-scale systems of equations with significantly reduced computational cost \cite{BGWSIREV}. In projection-based model reduction, the reduced systems are obtained via (Petrov-)Galerkin projection of the full-system equations onto low-dimensional---reduced---subspaces of the high-dimensional solution spaces corresponding to the full systems. If the large-scale systems contain nonlinear equations, then projection of the full-system equations onto reduced spaces typically is insufficient to obtain reduced systems that are computationally cheaper to solve than the full systems, because the nonlinear terms entail computations with costs that scale with the number of the degrees of freedom of the full system. The empirical interpolation method (\textsf{EIM}) \cite{barrault04-EIM,Maday2013,MADAY2015310}, and its discrete counter part, the discrete empirical interpolation method (\deim) \cite{DEIM,drmac-gugercin-DEIM-2016}, provide one solution to this problem by approximating the nonlinear terms of the nonlinear equations via sparse sampling. The nonlinear terms are evaluated at a few interpolation points---sampling points---and then all other components of the nonlinear terms are approximated via interpolation in low-dimensional subspaces. However, approximations via \ddeim have been shown to suffer from instabilities in certain situations, see, e.g., \cite{BeckettZhouMSThesis,GHAVAMIAN2017458,Mula}. Localization \cite{EftangLocal,peherstorfer13localized} and adaptation \cite{Peherstorfer15aDEIM,P18AADEIM} of the low-dimensional subspaces have been proposed as possible remedies. Another remedy that has been reported in the literature, and that typically is easier to implement in practice than localization and adaptation, is ``oversampling'' empirical interpolation so that the nonlinear terms are approximated via regression rather than via interpolation, which goes under the name of gappy proper orthogonal decomposition (\gpod) in the model reduction literature \cite{AstWWB08,CarlbergPhD2011,BeckettZhouMSThesis,Mula,ZimW16}. In this work, we consider the specific case where only noisy samples---observations---of the nonlinear terms are available and where \ddeim has been shown to be unstable, see, e.g., \cite{Mula}. We provide a probabilistic analysis that shows that \gpod with randomized samples leads to stable approximations in the presence of noise if more sampling points than basis vectors are used. 

Approximations based on regression, rather than interpolation, have been investigated in the context of model reduction. Missing point estimation (\textsf{MPE}) \cite{AstWWB04,AstWWB08} relies on \gpod \cite{Everson1995} to approximate nonlinear terms in model reduction. Several sampling point selection algorithms have been proposed for \textsf{MPE} and \gpod. The work \cite{WILLCOX2006208} formulates point selection as a sensor placement problem and proposes a greedy approach to find an approximate solution. Detailed analyses of point selection for \textsf{MPE}, and screening approaches to speedup point selection, are provided in \cite{AstWWB08}. The work by Zimmermann et al.~\cite{ZimW16} introduces a sampling strategy for \textsf{MPE} that is based on approximating eigenvalues for selecting sampling points and demonstrates that oversampling achieves higher accuracies in numerical experiments in computational fluid dynamics than \textsf{MPE} without oversampling. We will arrive at a special case of the approach presented in \cite{ZimW16} via perturbation bounds on eigenvalues introduced in \cite{IpsenBound}.  Carlberg et al.~\cite{carlberg2011efficient,Carlberg2013} introduce the Gauss-Newton with approximated tensors (\textsf{GNAT}) method that is based on Petrov-Galerkin projection and approximates the nonlinear terms via low-cost least-squares problems as in \gpod. The \textsf{GNAT} method and its performance based on regression has been investigated in the thesis \cite{CarlbergPhD2011}, where a greedy-based deterministic sampling strategy for selecting sampling points has been proposed. Zhou \cite{BeckettZhouMSThesis} introduces a deterministic sampling strategy for \gpod that exploits the dependency of the degrees of freedom of the full system to select sampling points. Regression via \gpod is then applied to multi-scale problems, where Zhou's sampling strategy with \gpod achieves lower errors than \deim via interpolation. The adaptive \deim (\textsf{ADEIM}), which adapts the \deim space from sparse samples of the nonlinear terms, is based on regression \cite{Peherstorfer15aDEIM,ZPW17SIMAXManifold,P18AADEIM}, even though regression is used for adaptation only and the nonlinear terms are approximated via interpolation once the \deim interpolants have been adapted. Other sampling strategies motivated by \deim and \gpod are investigated by Kutz et al.~\cite{2017arXiv170107569M}, who showed improvements for signal reconstruction \cite{PhysRevE.92.033304,MANOHAR2017162}. Greedy methods for sensor placement in the context of empirical interpolation are investigated in \cite{ARGAUD2018354,doi:10.1137/17M1157635}.

We consider \gpod in the specific setting where samples are polluted with noise. Noise is here to be understood in general terms, including perturbations that are typically modeled via random noise such as in turbulence, see, e.g., \cite{StochTurbulence}. It has been discussed in \cite{Mula} that the $L_2$ error of \ddeim approximations can grow with the dimension of the \ddeim space in presence of noise. The work \cite{Mula} proposes taking more sampling points than the dimension of the \ddeim space as a possible remedy and demonstrates on numerical results that this gives more stable results than \ddeim, i.e., that the error does not increase with the \ddeim dimension. We  build on the vast literature on \gpod and related methods \cite{AstWWB04,AstWWB08,ZimW16,carlberg2011efficient,Carlberg2013,BeckettZhouMSThesis,2017arXiv170107569M}. Our contribution is a probabilistic analysis that proves that in expectation with high probability \gpod with oversampling avoids the increase of the $L_2$ error with the dimension of the reduced space. For the analysis, we follow the work by Balzano et al.~\cite{5513344} and the work by Cohen et al.~\cite{Cohen2013} that provide approximation results for least-squares approximations, which we apply to \gpod with oversampling. Extensions to the work by Cohen et al. \cite{Cohen2013} have been introduced in \cite{MIGLIORATI2015167,refId0}. 
We then discuss a deterministic oversampling strategy and demonstrate with numerical results that a lower error with \gpod is achieved in the presence of noise compared to \ddeim that interpolates the nonlinear terms.

The structure of the paper is as follows. \Cref{sec:Prelim} briefly reviews \deim in the context of model reduction and numerically demonstrates on a toy example that \deim approximations are unstable if the nonlinear function evaluations are polluted with noise. \Cref{sec:Noise} and \Cref{sec:ODEIM} analyze \gpod with randomized samples and prove that oversampling avoids the stability issue in expectation with high probability. \Cref{sec:DetSampling} introduces a deterministic sampling strategy, which is then shown to achieve more accurate reduced models than \ddeim in \Cref{sec:NumRes}.

\section{Preliminaries and problem formulation}
\label{sec:Prelim}
This section briefly reviews \ddeim for approximating the nonlinear terms in reduced models and for recovering field data from few measurements and demonstrates, via an example, that \ddeim can become unstable in the presence of noise.

\subsection{Model reduction with empirical interpolation}
\label{sec:Prelim:MOR}
Consider a system of parametrized nonlinear equations
\begin{equation}
\bfA\bfx(\bfxi) + \bff(\bfx(\bfxi); \bfxi) = 0\,,
\label{eq:MOR:FOMMu}
\end{equation}
where $\bfx(\bfxi) \in \mathbb{R}^{\nh}$ is the state, $\bfxi \in \Dcal$ is a $\pdim$-dimensional parameter in the parameter domain $\Dcal$, $\bfA \in \mathbb{R}^{\nh \times \nh}$ is a constant matrix, and $\bff: \mathbb{R}^{\nh} \times \Dcal \to \mathbb{R}^{\nh}$ is a nonlinear function. Systems such as \cref{eq:MOR:FOMMu} typically arise after discretizing a PDE in the spatial domain, in which case  the matrix $\bfA$ corresponds to the linear operators of the underlying PDE and the nonlinear function $\bff$ to the nonlinear terms. In the following, we are interested in situations where the dimension $\nh \in \mathbb{N}$ of the state $\bfx(\bfxi)$ is large, which means that system \cref{eq:MOR:FOMMu} is potentially expensive to solve numerically, especially if these simulations need to be repeated for many parameter samples in outer-loop applications \cite{PWG17MultiSurvey} such as optimization, uncertainty quantification, and control.

A common approach to constructing a reduced model of the full system \cref{eq:MOR:FOMMu} is to use projection-based model reduction \cite{RozzaPateraSurvey,BGWSIREV}. Towards this goal, let the columns of the matrix $\bfX = [\bfx_1, \dots, \bfx_M] \in \mathbb{R}^{\nh \times M}$ be $M$ snapshots derived from the parameter samples $\bfxi_1, \dots, \bfxi_M \in \Dcal$ such that $\bfx_i = \bfx(\bfxi_i)$ for $i = 1, \dots, M$. Note that typically $M\leq N$. Further, let $\bfV = [\bfv_1, \dots, \bfv_{\nr}] \in \mathbb{R}^{\nh \times \nr}$ be an $\nr$-dimensional orthonormal basis constructed from the snapshot matrix $\bfX$. A common approach to obtaining $\bfV$ is to compute the singular value decomposition (\textsf{\textsf{SVD}}) of $\bfX$ and then to define $\bfV$ as the leading $r\leq M$ left singular vectors, as  done in proper orthogonal decomposition (\rpod). Then, the \rpod-Galerkin reduced model is obtained via projection
\begin{equation}
\stilde{\bfA}\tilde{\bfx}(\bfxi) + \bfV^T\bff(\bfV\tilde{\bfx}(\bfxi); \bfxi) = 0\,,
\label{eq:MOR:ROMMu}
\end{equation}
where $\stilde{\bfA} = \bfV^T\bfA\bfV$ is the reduced linear operator and $\tilde{\bfx}(\bfxi) \in \mathbb{R}^{\nr}$ is the reduced state.

Even though the  reduced state $\tilde{\bfx}(\bfxi)$ is in the $r$-dimensional subspace, evaluation of the reduced nonlinear term
$\bfV^T\bff(\bfV\tilde{\bfx}(\bfxi); \bfxi)$ in \cref{eq:MOR:ROMMu} still requires, first, lifting $\tilde{\bfx}(\bfxi)$  to the full dimension $\nh$, evaluating the original nonlinear term in this original dimension, and then projecting it down to the reduced dimension; thus, evaluating the reduced model
\cref{eq:MOR:ROMMu} still requires operations that scale with the dimension of the full model. This is called the lifting bottleneck in model reduction.

An effective remedy to the lifting bottleneck is the empirical interpolation method \cite{barrault04-EIM,DEIM}. The goal is to find an accurate approximation $\stilde{\bff}: \mathbb{R}^{\nd} \times \Dcal \to \mathbb{R}^{\nd}$ to $\bff$ that is computationally cheap to evaluate with cost independent of the dimension $\nh$. The empirical interpolation approximant $\stilde{\bff}$ has the form
\begin{equation}  \label{ftildedeim}
\stilde{\bff}(\tilde{\bfx}(\bfxi); \bfxi) = \bfU \bfc(\tilde{\bfx}(\bfxi); \bfxi)
\end{equation}
with $\tilde{\bfx}(\bfxi) \in \mathbb{R}^{\nd}$ and where $\bfc(\tilde{\bfx}(\bfxi); \bfxi) \in \mathbb{R}^\nd$ are the coefficients of the linear combination with the columns of $\bfU\in \mathbb{R}^{\nh \times \nd}$, which form a basis of an $\nd$-dimensional reduced space in which to approximate the function $\bff$ with $\nd \ll \nh$. \deim
achieves  the approximation \eqref{ftildedeim} by
 interpolating $\bff$ at selected components. Let
 $p_1, \dots, p_{\nd} \in \{1, \dots, \nh\}$ be the interpolation points (indices), i.e.,
 $ \bfe_{p_i}^T \bff(\bfV\tilde{\bfx}(\bfxi); \bfxi) = \bfe_{p_i}^T \stilde{\bff}(\tilde{\bfx}(\bfxi); \bfxi)$ for $i=1,2,\ldots,\nd$,
 where $\bfe_i \in \mathbb{R}^{\nh}$ denotes the $i$-th canonical unit vector. Let $\bfP = [\bfe_{p_1}, \dots, \bfe_{p_{\nd}}] \in \mathbb{R}^{\nh \times \nd}$ be the corresponding interpolation points (index selection) matrix. Then, the interpolation conditions are $\bfP^T \bff(\bfV\tilde{\bfx}(\bfxi); \bfxi) = \bfP^T \stilde{\bff}(\tilde{\bfx}(\bfxi); \bfxi)$, which, using \cref{ftildedeim}, lead to
 \begin{equation} \label{eq_deim}
 \stilde{\bff}(\tilde{\bfx}(\bfxi); \bfxi) = \bfU \bfc(\tilde{\bfx}(\bfxi); \bfxi) = \bfU (\bfP^T \bfU)^{-1}\bfP^T \bff(\bfV\tilde{\bfx}(\bfxi); \bfxi)\,,
 \end{equation}
 where $\bfc(\tilde{\bfx}(\bfxi); \bfxi) =  (\bfP^T \bfU)^{-1}\bfP^T \bff(\bfV\tilde{\bfx}(\bfxi); \bfxi)$.  In \cref{eq_deim},  $\stilde{\bff}$ is the \deim approximation of $\bff$.

The columns of $\bfU \in \mathbb{R}^{\nh \times \nd}$ are often taken as the \rpod basis of the
nonlinear snapshots  $\bff(\bfx(\bfxi_1); \bfxi_1), \dots, \bff(\bfx(\bfxi_M); \bfxi_M)$ with parameters $\bfxi_1, \dots, \bfxi_M \in \Dcal$. Note that $\bfU$ is orthonormal. The choice of the
selection operator $\bfP$ is motivated by the error bound
\begin{equation}\label{eq:DEIM-error0}
	\left\| \bff(\bfV\tilde{\bfx}(\bfxi); \bfxi) - \stilde{\bff}(\tilde{\bfx}(\bfxi); \bfxi) \right\|_2 \leq   \left\| (\bfP^T \bfU)^{-1} \right\|_2 \left\| (\Id-\bfU\bfU^T)\bff(\bfV\tilde{\bfx}(\bfxi); \bfxi)\right\|_2,
	\end{equation}
where $\| (\bfI-\bfU\bfU^T)\bff(\bfV\tilde{\bfx}(\bfxi); \bfxi)\|_2$ is the error due to the optimal approximation by orthogonal projection; see, \cite{barrault04-EIM,DEIM}. Therefore,
the selection operator $\bfP$ should choose indices such that $\| (\bfP^T \bfU)^{-1} \|_2$ is small. The \deim algorithm \cite{barrault04-EIM,DEIM} performs a greedy search to select the interpolation points.
The \qdeim point selection algorithm \cite{drmac-gugercin-DEIM-2016,drmac-saibaba-WDEIM-2018} based on the rank-revealing QR factorization is an alternative to this greedy-based point selection algorithms.
Combining the \deim approximation \cref{eq_deim} with the \rpod-Galerkin reduced model \cref{eq:MOR:ROMMu}, we obtain the \rpod-\deim-Galerkin reduced model
\begin{equation}
\stilde{\bfA}\tilde{\bfx}(\bfxi) + \bfV^T\bfU(\bfP^T\bfU)^{-1}\bfP^T\bff(\bfV\tilde{\bfx}(\bfxi); \bfxi) = 0\,,
\label{eq:MOR:ROMDEIM}
\end{equation}
where the $\nh - \nd$ components of $\bff(\bfV\tilde{\bfx}(\bfxi); \bfxi)$ that are different from the interpolation points $p_1, \dots, p_{\nd}$ are approximated via empirical interpolation. Thus, the  reduced model \cref{eq:MOR:ROMDEIM} requires evaluating the nonlinear function $\bff$ at only $\nd$ components, which typically leads to significant speedups compared to the \rpod-Galerkin reduced model \cref{eq:MOR:ROMMu} that requires evaluating the function $\bff$ at all $\nh$ components.

\begin{figure}
\centering
\begin{tabular}{cc}
{\LARGE\resizebox{0.45\columnwidth}{!}{\input{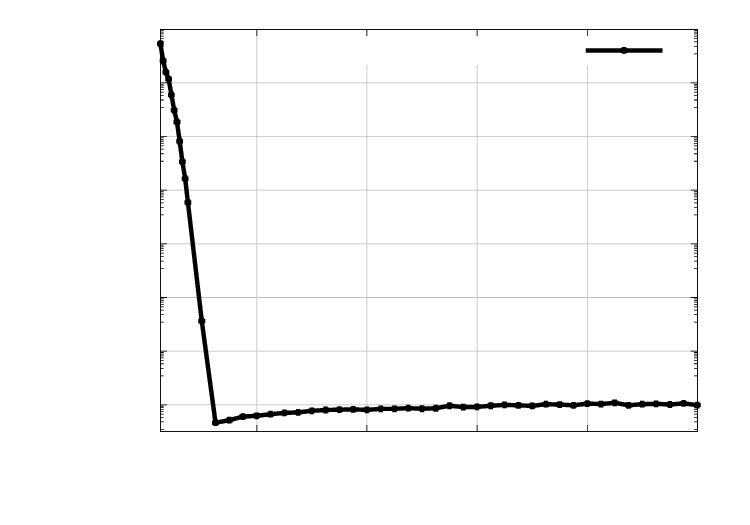}}} &
{\LARGE\resizebox{0.45\columnwidth}{!}{\input{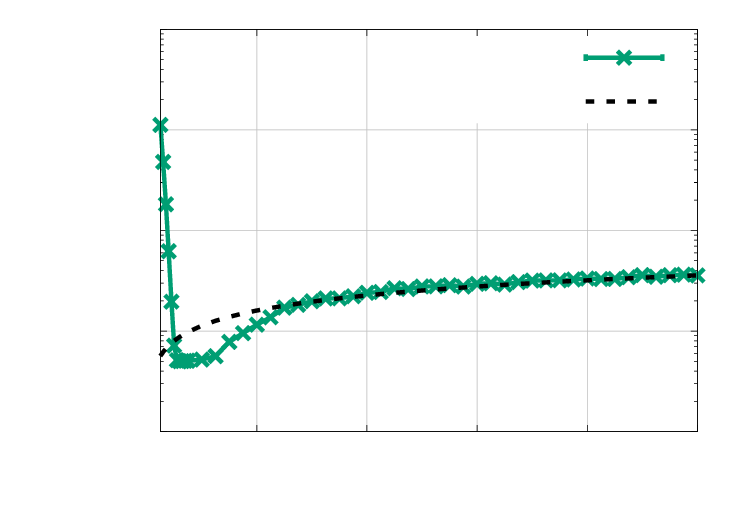}}}\\
\scriptsize (a) without noise & \scriptsize (b) with noise
\end{tabular}
\caption{The \ddeim is sensitive to noise in the sparse samples of the nonlinear function. In particular, the noise is amplified as the dimension $\nd$ of the reduced space is increased. A rate of $\sqrt{\nd}$ is numerically observed. Standard deviation of noise is $10^{-4}$.}
\label{fig:DEIMToy}
\end{figure}

\subsection{State field approximation from few measurements}
\label{sec:Prelim:Field}
Another use case of empirical interpolation and related methods, such as \gpod \cite{Everson1995}, is approximating state fields $\bfx: \Dcal \to \mathbb{R}^{\nh}$ from a few spatial measurements \cite{WILLCOX2006208,BuiDW,ARGAUD2018354,Mula}, where $\bfxi \in \Dcal$ is a parameter that defines the field $\bfx(\bfxi)$. Let $\bfU \in \mathbb{R}^{\nh \times \nd}$ be the reduced basis matrix constructed from snapshots $\bfx(\bfxi_1), \dots, \bfx(\bfxi_M)$ for $\bfxi_1, \dots, \bfxi_M \in \Dcal$ via, e.g., \rpod, and let $\bfP \in \mathbb{R}^{\nh \times \nd}$ be the interpolation points matrix derived from $\bfU$ with, e.g., the greedy algorithm \cite{barrault04-EIM,DEIM} and \qdeim \cite{drmac-gugercin-DEIM-2016,drmac-saibaba-WDEIM-2018}. Given are the measurements $\bfx_P(\bfxi) = \bfP^T\bfx(\bfxi)$ at the spatial coordinates corresponding to the components selected by $\bfP$ of a field $\bfx(\bfxi) \in \mathbb{R}^{\nh}$ with parameter $\bfxi \in \Dcal$. The field $\bfx(\bfxi)$ is unknown at all spatial coordinates except at the interpolation points given by $\bfP$. The \deim approximation of $\bfx(\bfxi)$ is then given by $\tilde{\bfx}(\bfxi) = \bfU(\bfP^T\bfU)^{-1}\bfx_P(\bfxi)$. For the ease of presentation, we follow the notation introduced in \Cref{sec:Prelim:MOR} for approximating nonlinear terms for model reduction; however, all what is presented in the following directly applies to state field approximation as well. We will revisit state field approximation in our numerical experiments in \Cref{sec:NumRes}.

\subsection{Instability of empirical interpolation in the presence of noise}
\label{sec:Prelim:Instability}
To approximate $\bff(\bfV\tilde{\bfx}(\bfxi); \bfxi)$ with \deim in the reduced model \cref{eq:MOR:ROMDEIM}, the function $\bff$ is evaluated (at least) at the components of $\bfV\tilde{\bfx}(\bfxi)$ corresponding to the interpolation points $p_1, \dots, p_{\nd}$, while all the other components are approximated via interpolation in the reduced space spanned by the columns of the basis matrix $\bfU$. We are interested in the situation where the function evaluations of $\bff$ at $\bfV\tilde{\bfx}(\bfxi)$ are noisy, in which case \deim approximations can become unstable, as demonstrated in, e.g., \cite{Mula}.

Consider the parametrized nonlinear function
\begin{equation}
\bff(\bfx; \xi) = \exp\left(-\frac{\left(\bfx - \xi\right)^2}{5 \times 10^{-3}}\right)\,,
\label{eq:Prelim:F}
\end{equation}
with the parameter $\xi \in \Dcal = [1, 3] \subset \mathbb{R}$. The components of $\bfx \in \mathbb{R}^{8192}$  are the equidistant points in $\Omega = [-2\pi, 2\pi]$. Note that all operations in \cref{eq:Prelim:F} are to be understood component-wise. Let $\xi_1, \dots, \xi_{2500}$ be the equidistant points in $\Dcal$ and let $\bff(\bfx; \xi_1), \dots, \bff(\bfx; \xi_{2500})$ be the nonlinear snapshots to derive a \deim interpolant $\stilde{\bff}$ of $\bff$ of dimension $\nd$ with the reduced basis matrix $\bfU$ and the \qdeim interpolation points matrix $\bfP$. We now approximate the function $\bff$ at the 2500 parameters $\xi_1^{\prime}, \dots, \xi_{2500}^{\prime} \in \Dcal$  uniformly sampled in the domain $\Dcal$. Note that the parameters $\xi_1^{\prime}, \dots, \xi_{2500}^{\prime}$ are different from the parameters $\xi_1, \dots, \xi_{2500}$ that were used to construct the reduced space and the interpolation points matrix. The \deim approximation
$\stilde{\bff}$
 of $\bff$ is
\[
\stilde{\bff}(\bfx; \xi_i^{\prime}) = \bfU(\bfP^T\bfU)^{-1}\bfP^T\bff(\bfx; \xi_i^{\prime})\,,
\]
for $i = 1, \dots, 2500$. The averaged relative state error
\begin{equation}
\frac{1}{2500}\sum_{i = 1}^{2500} \frac{\|\bff(\bfx; \xi_i^{\prime}) - \tilde{\bff}(\bfx; \xi_i^{\prime})\|_2}{\|\bff(\bfx; \xi_i^{\prime})\|_2}
\label{sec:Prelim:L2Error}
\end{equation}
versus the dimension $\nd$ of the \deim approximation is plotted in \Cref{fig:DEIMToy}a. The results indicate a fast decay of the \deim approximation error with the dimension $\nd$.

Let us now consider noisy evaluations of the function $\bff$. Therefore, let  $\bfepsilon$ be a random vector that has, as components, independent zero-mean Gaussian random variables with standard deviation $\sigma = 10^{-4}$. Define
\begin{equation}
\bff_{\bfepsilon}(\bfx; \xi) = \bff(\bfx; \xi) + \bfepsilon\,,
\label{eq:Prelim:NoisyF}
\end{equation}
so that the \deim approximation using the noisy function evaluations \cref{eq:Prelim:NoisyF} is
\[
\stilde{\bff}_{\bfepsilon}(\bfx; \xi) = \bfU(\bfP^T\bfU)^{-1}\bfP^T\bff_{\bfepsilon}(\bfx; \xi)\,.
\]
The plot in \Cref{fig:DEIMToy}b shows the averaged relative state error
\begin{equation}
\sum_{j = 1}^k \frac{1}{2500}\sum_{i = 1}^{2500} \frac{\|\bff(\bfx; \xi_i^{\prime}) - \stilde{\bff}_{\bfepsilon_j}(\bfx; \xi_i^{\prime})\|_2}{\|\bff(\bfx; \xi_i^{\prime})\|_2}\,,
\label{eq:Prelim:L2ErrorNoise}
\end{equation}
for $k = 10$ replicates of the \deim approximation $\tilde{\bff}_{\bfepsilon_j}$ that is derived from the noisy function evaluations \cref{eq:Prelim:NoisyF} with realization $\bfepsilon_j$ of the noise. The error bars indicate the minimum and maximum of the error over the replicates. Note that the error bars are barely visible, which means the variation over the replicates is small. The results indicate a stability issue of \deim in this case of noisy function evaluations because the error grows with the dimension $\nd$ of the reduced space. The result illustrates an error growth with a rate $\sqrt{\nd}$ with the dimension $\nd$. Similar observations are made in \cite{Mula}.

\begin{remark}
The term ``instability'' has various meanings in numerical analysis. In the following, the term ``instability'' refers to the specific phenomenon that the \deim approximation error $\|\bff(\bfx; \bfxi) - \stilde{\bff}_{\bfepsilon}(\bfx; \bfxi)\|_2$ in the Euclidean norm grows with the dimension $\nd$ of the reduced space if noisy function evaluations \eqref{eq:Prelim:NoisyF} (or noisy measurements in the context of state field approximation in \Cref{sec:Prelim:Field}) are used; see Figure~\ref{fig:DEIMToy}b.
\end{remark}

\section{Amplification of noise in \deim}
\label{sec:Noise}
We provide an upper bound on the amplification of the noise in \deim approximations, and a theoretical explanation of the numerical observation in \Cref{fig:DEIMToy}b. The bound \cref{lemma1-claim3}, which we prove in the following, shows that the error cannot increase faster than with rate $\sqrt{\nd}$, which is the rate observed in  \Cref{fig:DEIMToy}b. We also provide a formula for the expected value of the \deim error vector and reveal the structure of the error ellipsoid. A similar bound as \cref{lemma1-claim3} has been presented in \cite{Mula}.

To simplify the exposition, we drop the dependence on the state $\bfx$ and the parameter $\bfxi$ of $\bff(\bfx; \bfxi)$,  and abbreviate it as $\bff(\bfx; \bfxi) = \bff$.  Similarly, the \deim approximant will be abbreviated as $\stilde{\bff}$. The noisy counterparts of $\bff$ and $\stilde{\bff}$ are $\bff_{\bfepsilon} = \bff + \bfepsilon$ and $\stilde{\bff}_{\bfepsilon}$, respectively, where $\bfepsilon$ is a zero-mean Gaussian vector with independent components with standard deviation $\bfsigma = [\sigma_1, \dots, \sigma_{\nh}]^T$. 
\begin{lem}\label{lm:InterBound}
Define the error of the DEIM approximation $\tilde{\bff}_{\bfepsilon}$ from noisy function evaluations as $\bfr_{\bfepsilon} = \bff - \tilde{\bff}_{\bfepsilon} = \bff - \bfU(\bfP^T\bfU)^{-1}\bfP^T\bff_{\bfepsilon}$ and the error of the approximation $\tilde{\bff}$ with noise-free function evaluations as $\bfr = \bff - \bfU(\bfP^T\bfU)^{-1}\bfP^T\bff$. Then, the expected value of the error $\bfr_{\bfepsilon}$ corresponding to noisy function evaluations equals $\bfr$, i.e., $\mathbb{E}_{\bfepsilon}[\bfr_{\bfepsilon}] = \bfr$, where the expectation is taken over the noise. The standard deviation of $\bfr_{\bfepsilon}$ satisfies
\begin{equation}\label{lemma1-claim2}
\mathbb{E}_{\bfepsilon}\left[\| \bfr_{\bfepsilon} - \mathbb{E}_{\bfepsilon}[\bfr_{\bfepsilon}]\|_2\right] \leq \sqrt{\mathbb{E}_{\bfepsilon}\left[\| \bfr_{\bfepsilon} - \mathbb{E}_{\bfepsilon}[\bfr_{\bfepsilon}]\|_2^2\right] }  \leq \sqrt{n} \| (\bfP^T\bfU)^{-1}\|_2 \|\bfP^T\bfsigma\|_{\infty}\,.
\end{equation}
Thus, the error is bounded in expectation as
\begin{equation}\label{lemma1-claim3}
\mathbb{E}_{\bfepsilon}\!\left[\|\bff - \bfU(\bfP^T\bfU)^{-1}\bfP^T\bff_{\bfepsilon}\|_2\right] \leq \!
\| (\bfP^T\bfU)^{-1}\|_2\!  \left(\|\bff - \bfU\bfU^T\bff\|_2 + \sqrt{\nd}\|\bfP^T\bfsigma\|_{\infty}\right).
\end{equation}
\end{lem}
\begin{proof}
Using the linearity of the expectation, the error formula for the \deim projection, and  the assumptions on the noise, namely, $\mathbb{E}_{\bfepsilon}\left[ \bfepsilon \right] = 0$, we obtain
\begin{align}
 \mathbb{E}_{\bfepsilon}\left[ \bff - \bfU(\bfP^T\bfU)^{-1}\bfP^T\bff_{\bfepsilon}\right] &=
(\bfI - \bfU(\bfP^T\bfU)^{-1}\bfP^T)\bff - \bfU(\bfP^T\bfU)^{-1} \bfP^T \mathbb{E}_{\bfepsilon}\left[ \bfepsilon \right]\\
 &= (\bfI - \bfU(\bfP^T\bfU)^{-1}\bfP^T)\bff,
\nonumber
\end{align}
which establishes $\mathbb{E}_{\bfepsilon}\left[\bfr_{\bfepsilon}\right] = \bfr$ as claimed.
The norm of $\|\mathbb{E}_{\bfepsilon}[\bfr_{\bfepsilon}]\|_2$ is bounded as
\begin{equation}
\left\|\mathbb{E}_{\bfepsilon}\left[ \bff - \bfU(\bfP^T\bfU)^{-1}\bfP^T\bff_{\bfepsilon}\right] \right\|_2 \leq
\left\|(\bfP^T\bfU)^{-1}\right\|_2 \left\|\bff - \bfU\bfU^T\bff\right\|_2,
\end{equation}
which is the same upper bound as in \eqref{eq:DEIM-error0} for the noise-free case.
The covariance matrix of the error $\bfr_{\bfepsilon}$ is
\begin{eqnarray}
\bfC = \mathbb{E}_{\bfepsilon}\left[ \bfU(\bfP^T\bfU)^{-1}\bfP^T{\bfepsilon}  {\bfepsilon}^T \bfP (\bfP^T\bfU)^{-T}\bfU^T\right]   &=&
\bfU(\bfP^T\bfU)^{-1}\bfP^T\Sigma^2 \bfP (\bfP^T\bfU)^{-T} \bfU^T  \nonumber \\ &= &
\bfU(\bfP^T\bfU)^{-1}\Sigma_{\bfP}^2 (\bfP^T\bfU)^{-T} \bfU^T\,,\label{eq:SigmaP}
\end{eqnarray}
where $\Sigma_{\bfP}^2 =\bfP^T\Sigma^2 \bfP = \mathrm{diag}(\sigma_{p_i}^2)_{i=1}^n$ with $\sigma_{p_i}$ being the standard deviation of the $p_i$-th component of $\bfepsilon$. The covariance $\bfC $ is positive semidefinite. Its nonzero eigenvalues $\lambda_i^2$ (that correspond to the invariant space spanned by $\bfU$) can be enumerated so that
$$
\lambda_i^2 = \sigma_{p_i}^2 \vartheta_i^2,\;\; \mbox{where} \;\;\frac{1}{\|\bfP^T\bfU\|_2} \leq \vartheta_i \leq \| (\bfP^T\bfU)^{-1}\|_2,\;\;i=1,\ldots, n .
$$
This is an application of the Ostrowski theorem \cite[Theorem 4.5.9]{Horn:2012:MA:2422911}; it identifies the bounds of the amplification factors $\vartheta_i$'s of the corresponding standard deviations.
The spectral structure of $\bfC$ (and thus the error ellipsoid) can be explicitly revealed using the \textsf{SVD} $\bfSigma_{\bfP}^{-1}(\bfP^T\bfU)=\bfPhi \bfOmega\bfPsi^T$
($\bfPhi $,  $\bfPsi$ orthogonal matrices of singular vectors, $\bfOmega$ diagonal matrix of singular values), which yields
$\bfC = (\bfU \bfPsi) \bfOmega^{-2} (\bfU \bfPsi)^T$.

Next, of interest is the variance $\mathbb{E}_{\bfepsilon}[\| \bfr_{\bfepsilon} - \mathbb{E}_{\bfepsilon}[\bfr_{\bfepsilon}]\|_2^2]$ of $\bfr_{\bfepsilon}$, for which follows that
$$
\mathbb{E}_{\bfepsilon}\left[\| \bfr_{\bfepsilon} - \mathbb{E}_{\bfepsilon}[\bfr_{\bfepsilon}]\|_2^2\right] = \mathrm{Trace}(\bfC)=\sum_{i=1}^n \sigma_{p_i}^2 \vartheta_{i}^2 \leq n \| (\bfP^T\bfU)^{-1}\|_2^2 \max_{i}\sigma_{p_i}^2\,,
$$
which shows $\mathbb{E}_{\bfepsilon}\left[\| \bfr_{\bfepsilon} - \mathbb{E}_{\bfepsilon}[\bfr_{\bfepsilon}]\|_2^2\right] \leq n \|\bfP^T\bfsigma\|_{\infty}^2 \|(\bfP^T\bfU)^{-1}\|_2^2$. In addition, by deploying Jensen's inequality and taking square root, we obtain \cref{lemma1-claim2}.
Then, combining the above estimates and the triangle inequality yields
\begin{eqnarray*}
\mathbb{E}_{\bfepsilon}\left[\| \bfr_{\bfepsilon} \|_2\right] &\leq& \mathbb{E}_{\bfepsilon}\left[\| \bfr_{\bfepsilon} - \mathbb{E}_{\bfepsilon}[\bfr_{\bfepsilon}]\|_2\right]  +\|\mathbb{E}_{\bfepsilon}\left[\bfr_{\bfepsilon} \right]\|_2
\\ &\leq&  \| (\bfP^T\bfU)^{-1}\|_2 ( \|\bff - \bfU\bfU^T\bff\|_2 + \sqrt{n} \|\bfP^T\bfsigma\|_{\infty}).
\end{eqnarray*}
 which proves \cref{lemma1-claim3}.
\end{proof}
\begin{remark}
\cool{\Cref{lm:InterBound} assumes that the noise
 $\bfepsilon$ is a zero-mean Gaussian vector with independent components. If   $\bfepsilon$ does not have independent components, i.e., if the noise error covariance matrix $\bfSigma^2$ is not diagonal,  we can write the spectral decomposition of its submatrix $\bfSigma_{P}^2
  = \bfP^T\bfSigma^2\bfP$
 as $\bfSigma_P^2 = W_{\bfP} \bfD_{\bfP}^2 W_{\bfP}^T$
 where $\bfD_{\bfP}=\mathrm{diag}(\bfd_{i})_{i=1}^n$.
Using this spectral decomposition, we replace
\eqref{eq:SigmaP} with
$$
\bfC =\bfU(\widetilde{\bfP}^T\bfU)^{-1}\bfD_{\bfP}^2 (\widetilde{\bfP}^T\bfU)^{-T} \bfU^T
$$
where $\widetilde{\bfP} = \bfP \bfW_{\bfP}$. Since $\bfW_{\bfP}$ is an orthogonal matrix, $\|\widetilde{\bfP}^T\bfU\|_2=\|\bfP^T\bfU\|_2$ and $\|(\widetilde{\bfP}^T\bfU)^{-1}\|_2=\|(\bfP^T\bfU)^{-1}\|_2$. Then, the rest of the proof of \Cref{lm:InterBound} follows as before by replacing $\bfP$
with
$\widetilde{\bfP}$  and  $\sigma_{p_i}$ with $\bfd_{i}$, for $i=1,\ldots, n$. Finally, the upper bounds \eqref{lemma1-claim2} and
\eqref{lemma1-claim3} hold true by replacing $\bfP^T\bfsigma$  with the vector $\left[~\bfd_1,\ldots, \bfd_n\right]^T$, where now  $\bfd_i^2$'s are the variances along the principal components (eigenvectors of $\bfSigma_{\bfP}$). Moreover, in both \eqref{lemma1-claim2} and
\eqref{lemma1-claim3},
the term $\sqrt{n}\|\bfP^T \bfsigma\|_\infty$  can be replaced with $\|\bfSigma_{\bfP}\|_F =
\sqrt{\bfd_1^2+\cdots+\bfd_n^2}
$.}
\end{remark}
\begin{cor}
	If  the selection operator $\bfP$ is based on  the quasi-optimal point selection introduced in \cite[Lemma~2.1]{drmac-saibaba-WDEIM-2018},
then \cref{lemma1-claim2} and \cref{lemma1-claim3} hold with the bound
	\begin{equation}
	\|(\bfP^T\bfU)^{-1}\|_2 \leq \sqrt{1 + \eta^2 n(N - n)}\,.
	\label{eq:DERM:NoiseLemma:Inter3}
	\end{equation}
	where $\eta\geq 1$ is a tuning parameter.
\end{cor}

\begin{remark}

	Note how in \cref{lemma1-claim3}, with increasing column dimension $n$ of the matrix $\bfU$, the POD projection error $\|\bff - \bfU\bfU^T\bff\|_2$ monotonically decreases toward zero and, at the same time, the norm of the sampling operator $\|(\bfP^T\bfU)^{-1}\|_2$ approaches one, while the contribution of the noise grows as $\sqrt{n} \|\bfP^T\bfsigma\|_{\infty}$, taking over the leading term. The effect of the noise dominating the error is seen in \Cref{fig:DEIMToy}.

\end{remark}
\begin{remark}

From \Cref{lm:InterBound}, it follows  that it is desirable that a \deim selection operator  avoids components of $\bff$ with noise with high variance. Such a strategy may help slow the noise buildup.  If we denote by $\mathcal{J}$ undesirable indices and set $\mathcal{J}^c=\{1,\ldots, \nh\}\setminus\mathcal{J}$, then we can run the \qdeim selection on the submatrix $\bfU(\mathcal{J}^c,:)$; for details we refer the reader to  \cite[\S 3.]{drmac-gugercin-DEIM-2016}.

\end{remark}

\section{Stability of \gpod with randomized samples}
\label{sec:ODEIM}
Given the reduced basis matrix $\bfU \in \mathbb{R}^{\nh \times \nd}$, the \deim selects $\nd$ interpolation points, i.e., $\bfP^T\bfU$ is a square matrix. In this section, we investigate oversampling in the sense that more sampling points $\np > \nd$ than the dimension $\nd$ of the reduced space spanned by the columns of $\bfU$ are used. Taking more sampling points than the dimension of the space goes under the name of gappy proper orthogonal decomposition (\gpod), which has been introduced in \cite{Everson1995} and is used in the context for model reduction in \cite{AstWWB04,AstWWB08}. We now show that with \gpod, the noise amplification that was observed in \Cref{sec:Noise} can be avoided in expectation with high probability if sampling points are randomized and if more sampling points than basis vectors are used.

\subsection{\gpod}
Consider $p_1, \dots, p_{\np} \in \{1, \dots, \nh\}$,  pairwise distinct sampling points  with $\np > \nd$, i.e., the number of sampling points $\np$ is larger than the dimension $\nd$ of the space spanned by the columns of the basis matrix $\bfU$. Then,  the corresponding \gpod approximation of $\bff$ is
\[
\shat{\bff} = \bfU(\bfP^T\bfU)^\dagger\bfP^T\bff\,,
\]
where ${\bfM}^\dagger$ denotes the Moore-Penrose inverse of  $\bfM$, i.e.,
$
\bfM^\dagger =\left(\bfM^T\bfM\right)^{-1}\bfM^T$,
assuming $\bfM$ has linearly independent columns.
In contrast to the \deim approximation $\stilde{\bff}$ in \cref{eq_deim},
 the \gpod approximation $\shat{\bff}$ is obtained via regression and therefore does not necessarily interpolate $\bff$ at the sampling points $p_1, \dots, p_{\np}$. In the case of noise-free sampling, the error of the \gpod approximation in the Euclidean norm satisfies (see, e.g., \cite[Proposition~2.1]{ZimW16})
\begin{equation}
\|\bff - \shat{\bff}\|_2 \leq \|(\bfP^T\bfU)^\dagger\|_2 \|\bff - \bfU\bfU^T\bff\|_2\,,
\label{eq:DERMBound}
\end{equation}
where $\|(\bfP^T\bfU)^\dagger\|_2$ quantifies the effect of the sampling points and $\|\bff - \bfU\bfU^T\bff\|_2$ relates to the approximation quality of the space spanned by $\bfU$; cf.~the \deim error bound \cref{eq:DEIM-error0}.

\subsection{Probabilistic analysis of \gpod}
\label{sec:ODEIM:ProbAnalysis}
We now investigate the error of the \gpod approximation $\shat{\bff}$ when the sampling points $p_1, \dots, p_{\np}$ are selected uniformly with replacement from $\{1, \dots, \nh\}$. Note that the following analysis is developed for uniform sampling \emph{with replacement} as in the work by Balzano et al.~\cite{5513344}. Parts of the following analysis are an application of the work by Cohen et al.~\cite{Cohen2013}.

To set up the analysis, we define the coherence of a subspace $\Ucal = \operatorname{span}(\bfU)$ as
\[
\mu(\Ucal) = \frac{\nh}{\nd}\max_{i = 1, \dots, \nh} \| \rowu_i^T\|_2^2\,,
\]
where the columns of $\bfU \in \mathbb{R}^{\nh \times \nd}$ form an orthonormal basis for $\Ucal$ and $\rowu_i^T$ is the $i$-th row of $\bfU$, see, e.g., \cite[Definition~1.2]{Candes2009}. Intuitively speaking, coherence measures if there are certain coordinate directions that carry significantly more information than other directions. Note that $\max_{i = 1, \dots, \nh} \| \rowu_i^T\|_2^2\geq \nd/\nh$. The following result  from \cite[Lemma~3]{5513344} will be used in our analysis.
\begin{lem}
Let the points $p_1, \dots, p_{\np}$ be uniformly sampled from $\{1, \dots, \nh\}$ with replacement and let $\bfP$ be the corresponding sampling points matrix. Moreover, let  $\delta \in (0,1]$ such that $\np \geq (8/3) \nd \mu(\Ucal)\log(2\nd/\delta)$  and set $\gamma = \sqrt{\frac{8\nd\mu(\Ucal)}{3\np}\log(2\nd/\delta)}$. Then
\[
\left\|\left(\left(\bfP^T\bfU\right)^T\left(\bfP^T\bfU\right)\right)^{-1}\right\|_2 \leq \frac{\nh}{(1 - \gamma)\np}
\]
with probability at least $1 - \delta$.
\label{lm:EigenBound}
\end{lem}
The following lemma states that $\|(\bfP^T\bfU)^\dagger\|_2$ can be bounded by a constant, arbitrarily close to $1$,  with high probability if a sufficiently large number of sampling points is used, which means that \gpod with randomized samples is well-posed with high probability.

\begin{lem}
Consider the same setup as in \Cref{lm:EigenBound} and set
$\widehat{\gamma}=\sqrt{m}\gamma$. If $m$ is such that, for $K\geq 1$,
\begin{equation}\label{sqrt-m-lover}
\sqrt{m} \geq \frac{1}{2} \widehat{\gamma} + \frac{1}{2}\sqrt{\widehat{\gamma}^2 + \frac{4N}{K^2}},
\end{equation}
then
\[
\|(\bfP^T\bfU)^\dagger\|_2 \leq {K},
\]
with probability at least $1 - \delta$.
\end{lem}
\begin{proof}
{Since
	\[
	\|(\bfP^T\bfU)^\dagger\|_2  = \frac{1}{\sigma_{\text{min}}(\bfP^T\bfU)} = \sqrt{\frac{1}{\lambda_{\text{min}}((\bfP^T \bfU)^T(\bfP^T\bfU))}} =  \sqrt{\|((\bfP^T\bfU)^T(\bfP^T\bfU))^{-1}\|_2}
	\]
holds, \Cref{lm:EigenBound} yields that with probability at least $1 - \delta$
	\begin{equation} \label{eq:ptupinv}
	\|(\bfP^T\bfU)^\dagger\|_2  \leq \sqrt{\frac{N}{(1-\gamma)m}}.
	\end{equation}
Then, the task is to choose $m$ so that $N/\big((1-\gamma)m\big) \leq K^2$. To that end, set $\widehat{\gamma}=\sqrt{m}\gamma$. By the assumption of \Cref{lm:EigenBound}, $\gamma \leq 1$ and thus $\sqrt{m} \geq \widehat{\gamma}$.  The desired inequality becomes
$$
\nh \leq K^2 ( m - \sqrt{m}\widehat{\gamma}),\;\;\mbox{i.e.,}\;\;K^2 x^2 - K^2\widehat{\gamma}x-N \geq 0,\;\;\mbox{where}\;\;x=\sqrt{m}\geq\widehat{\gamma}.
$$
The smaller root of the above parabola is negative and the larger one, then, provides the  desired lower bound \cref{sqrt-m-lover}.
}
\end{proof}

The following bound will be helpful in establishing the main result.
\begin{lem}\label{A:LEM:4.3}
Consider the same setup as in \Cref{lm:EigenBound}.  Let  {$\bfg \in \mathbb{R}^{\nh}$}, and let $\alpha$ be the acute angle between $\bfg$ and the range of $\bfU$. Then,
{
\begin{equation}
\mathbb{E}_{P}\left[\big\|(\bfP^T\bfU)^\dagger\bfP^T\bfg\big\|_2\right] \leq  \min\left\{  {\frac{1}{\sqrt{1 - \gamma}}}, \frac{\sqrt{\cos^2\alpha + \frac{\nd}{\np}\mu(\Ucal)}}{1-\gamma}\right\} \|\bfg\|_2
\label{A:eq:lmPTU}
\end{equation}
}with probability at least $1 - \delta$,  where the expected value $\mathbb{E}_{P}$ is with respect to the uniform distribution of the sampling points.
\end{lem}
\begin{proof}
We first apply submultiplicativity to obtain
\begin{equation}  \label{ub1}
\|(\bfP^T\bfU)^\dagger\bfP^T\bfg\|_2 \leq \|(\bfP^T\bfU)^\dagger\|_2 \|\bfP^T\bfg\|_2\,.
\end{equation}
Using  \cref{eq:ptupinv} and \Cref{lm:EigenBound}, we have, with probability at least $1 - \delta$, that
\begin{equation}
\mathbb{E}_{P}\left[ \big\|(\bfP^T\bfU)^\dagger \bfP^T\bfg\big\|_2 \right]
\leq \sqrt{\frac{N}{(1-\gamma)m}} \,\,\mathbb{E}_{P}\left[\big\|\bfP^T\bfg\big\|_2\right]\,.
\label{A:eq:DERM:LemmaPseudo:Inter2}
\end{equation}
Consider now the expected value $\mathbb{E}_{P}[\|\bfP^T\bfg\|^2_2]$ and note that we use the squared Euclidean norm $\|\cdot\|^2_2$.
Let $g_j$ denote the $j$-th component of $\bfg$ for $j = 1, \dots, \nh$. Also let
 $\mathbb{I}_{p_i=j}$ denote the indicator function that is $1$ if $p_i=j$ and 0 otherwise.
 Note that the probability that $p_i = j$ is  $1/\nh$ because a uniform distribution with replacement is used for selecting the sampling points, and thus
 $\mathbb{E}[\mathbb{I}_{p_i=j}] = 1/\nh$.
Then,
\begin{equation}
\mathbb{E}_{P}\left[\big\|\bfP^T\bfg\big\|_2^2\right] = \mathbb{E}_{P}\left[  \sum_{i = 1}^{\np} \sum_{j = 1}^{\nh} g_j^2\,\mathbb{I}_{p_i=j}\right]\,
\label{A:eq:DERM:LemmaPseudo:Inter1a}
= \frac{m}{\nh} \|\bfg\|_2^2\,.
\end{equation}
Applying Jensen's inequality to \cref{A:eq:DERM:LemmaPseudo:Inter1a} yields
$
\mathbb{E}_{P}[\|\bfP^T\bfg\|_2] \leq \sqrt{\frac{m}{N}}\|\bfg\|_2\,,
$
which,  combined with \cref{A:eq:DERM:LemmaPseudo:Inter2} implies
\begin{equation}
\mathbb{E}_{P}\left[\big\|(\bfP^T\bfU)^\dagger\bfP^T\bfg\big\|_2\right] \leq  {\frac{1}{\sqrt{1 - \gamma}}} \| \bfg \|_2  \label{A:eq:min1},
\end{equation}
proving the upper bound in  \cref{A:eq:lmPTU} for the first input of the $\min$ function.
To prove \cref{A:eq:lmPTU} for the second input of the $\min$ function, we first apply submultiplicativity to obtain
\begin{equation}  \label{ub2}
\|(\bfP^T\bfU)^\dagger\bfP^T\bfg\|_2 \leq \left\|\left((\bfP^T\bfU)^T(\bfP^T\bfU)\right)^{-1}\right\|_2 \|(\bfP^T\bfU)^T\bfP^T\bfg\|_2\,.
\end{equation}
With \Cref{lm:EigenBound}, we have, with probability at least $1 - \delta$, that
\begin{multline}
\mathbb{E}_{P}\left[\left\|\left((\bfP^T\bfU)^T(\bfP^T\bfU)\right)^{-1}\right\|_2 \big\|(\bfP^T\bfU)^T\bfP^T\bfg\big\|_2\right] \\
\leq \frac{N}{(1-\gamma)m} \mathbb{E}_{P}\left[ \big\|(\bfP^T\bfU)^T\bfP^T\bfg\big\|_2\right]\,.
\label{A:eq:DERM:LemmaPseudo:Inter2a}
\end{multline}
Let $\bla \bfv\,, \bfw \bra_2 =  \bfv^T\bfw$ denote the Euclidean inner product and  consider  the expected value $\mathbb{E}_{P}[\|(\bfP^T\bfU)^T\bfP^T\bfg\|^2_2]$. Using  the  linearity of $\bla\cdot,\cdot\bra_2$ and $\mathbb{E}_{P}\left[\cdot\right]$, we obtain
\begin{eqnarray}
&&\mathbb{E}_{P}[\langle(\bfP^T\bfU)^T\bfP^T\bfg,(\bfP^T\bfU)^T\bfP^T\bfg\rangle_2]=
\mathbb{E}_{P}\left[\sum_{j=1}^m \sum_{k=1}^m \bla g_{p_j} \rowu_{p_j}^T,g_{p_k} \rowu_{p_k}^T \bra_2\right]  \nonumber
\\
&&~\quad = \mathbb{E}_{P}\left[\sum_{j=1}^m \sum_{k=1}^m \bla \sum_{\ell=1}^\nh g_{\ell} \rowu_{\ell}^T \mathbb{I}_{p_j=\ell} ,  \sum_{s=1}^\nh g_{s} \rowu_{s}^T \mathbb{I}_{p_k=s} \bra_2\right]  \nonumber \\
&&~\quad =
\sum_{j=1}^m \sum_{k=1}^m  \sum_{\ell=1}^\nh \sum_{s=1}^\nh g_{\ell}g_{s} \bla\rowu_{\ell}^T  ,    \rowu_{s}^T\bra_2\, \mathbb{E}_{P}\left[\mathbb{I}_{p_j=\ell}  \mathbb{I}_{p_k=s} \right] . \label{A:eq:E_PUG}
\end{eqnarray}
 Note that, for $k\neq j$, by independence of the $j$th and the $k$th drawing with replacement
\begin{equation}
\mathbb{E}_P\!\left[\mathbb{I}_{p_j=\ell}\mathbb{I}_{p_k=s}     \right] =
  \frac{1}{\nh^2}\label{A:eq:E:P:split}\,.
\end{equation}
Now we can use \cref{A:eq:E:P:split} to split \eqref{A:eq:E_PUG}.  For $k\neq j$, we obtain
$$
\sum_{j=1}^m \sum_{\substack{k=1\\k\neq j}}^m  \sum_{\ell=1}^\nh \sum_{s=1}^\nh g_{\ell}g_{s} \bla\rowu_{\ell}^T  ,    \rowu_{s}^T\bra_2 \frac{1}{\nh^2}=
\frac{m^2-m}{N^2}\bla\bfU \bfU^T \bfg,\bfg\bra_2= \frac{m^2-m}{N^2} \|\bfU \bfU^T \bfg\|_2^2.
$$
The remaining terms with $j=k$ contribute to \eqref{A:eq:E_PUG} with
\begin{align*}
\sum_{j=1}^m\sum_{\ell=1}^\nh \sum_{s=1}^\nh g_{\ell}g_{s} \bla\rowu_{\ell}^T  ,    \rowu_{s}^T\bra_2\, \mathbb{E}_{P}\left[\mathbb{I}_{p_j=\ell} \,\mathbb{I}_{p_j=s} \right]
& = \frac{1}{N}\left[\sum_{j=1}^m  \sum_{\ell=1}^\nh g_{{\ell}}^2 \|\rowu_{{\ell}}^T\|_2^2 \,
\right] \\&= \frac{\np}{\nh}  \sum_{\ell=1}^\nh g_{{\ell}}^2 \|\rowu_{{\ell}}^T\|_2^2\,,
\end{align*}
where we use $\mathbb{E}_P[\mathbb{I}_{p_j = \ell}\mathbb{I}_{p_j = s}] = 1/N$ for $s = \ell$ and 0 otherwise.
Altogether, we have
\begin{eqnarray}
&&\nonumber\mathbb{E}_{P}[\|(\bfP^T\bfU)^T\bfP^T\bfg\|^2_2] = \frac{\np^2-\np}{N^2}\|\bfU \bfU^T \bfg\|_2^2 + \frac{\np}{\nh} \| \mathrm{diag}(\bfg) \bfU\|_F^2  \\
&&~\qquad \leq \nonumber
\frac{\np^2-\np}{\nh^2} {\|\bfU \bfU^T \bfg\|_2^2} + \frac{\np}{\nh} \frac{n}{\nh} \mu(\Ucal) \|\bfg\|_2^2
\;\;\qquad (\mbox{here}\; \frac{n}{\nh} \mu(\Ucal)=\max_{\ell=1,\ldots,\nh}\|\rowu_{\ell}^T\|_2^2) \\
&& ~\qquad= \frac{\np^2-\np}{\nh^2}  \|\bfg\|_2^2 {\cos^2\alpha} + \frac{\np}{\nh} \frac{\nd}{\nh} \mu(\Ucal) \|\bfg\|_2^2\,,
\label{eq:all}
\end{eqnarray}
where $\cos^2\alpha=\|\bfU^T \bfg\|_2^2/\|\bfg\|_2^2$ and $\operatorname{diag}(\bfg) \in \mathbb{R}^{\nh \times \nh}$ is the diagonal matrix with the components of $\bfg$ on its diagonal. Hence, applying Jensen's inequality to \eqref{eq:all} and then using \cref{A:eq:DERM:LemmaPseudo:Inter2a} yields
$$
\mathbb{E}_{P}[\|(\bfP^T\bfU)^{\dagger}\bfP^T\bfg\|_2] \leq \frac{\|\bfg\|_2}{1-\gamma} \sqrt{{\cos^2\!\alpha}\! +\! \frac{\nd}{\np} \mu(\Ucal)} =
\frac{\|\bfg\|_2}{1-\gamma} \sqrt{{\cos^2\!\alpha} + \frac{\nh}{\np} \max_{\ell=1,\ldots,\nh}\|\rowu_{\ell}^T\|_2^2} .
$$
Combining this final inequality with \cref{A:eq:min1}  yields the desired result \eqref{A:eq:lmPTU}.
\end{proof}
\begin{remark}
The $\min$ function in the upper bound \cref{A:eq:lmPTU} results from using two different upper bounds
for $\|(\bfP^T\bfU)^\dagger\bfP^T\bfg\|_2$, one as in \cref{ub1} and the other one as in \cref{ub2}. The latter is employed by \cite{5513344} for the special case where $\bfg$ is orthogonal to the range of $\bfU$. A technical calculation shows that
while  the first input in the $\min$ function is expected to be the smaller of the two for small $\np$, the second input is expected to be the smaller one  for large $\np$.
\end{remark}
\begin{remark}\label{A:rem:gpU}
	{
If $\bfg$ is orthogonal to the range of $\bfU$,  in \cref{A:eq:lmPTU}  the term $\cos \alpha=0$   and
the upper bound simplifies to
$$
\mathbb{E}_{P}\left[\big\|(\bfP^T\bfU)^\dagger\bfP^T\bfg\big\|_2\right] \leq  \min\left\{  {\frac{1}{\sqrt{1 - \gamma}}}, \frac{\sqrt{ \frac{\nd}{\np}\mu(\Ucal)}}{1-\gamma}\right\} \|\bfg\|_2.
$$
Thus, \Cref{A:LEM:4.3} contains \cite[Lemma 2]{5513344} as a special case.
Indeed, in this special case, the expectation of $(\bfP^T\bfU)^T\bfP^T\bfg$ is zero because
$$
\mathbb{E}_{P}[(\bfP^T\bfU)^T\bfP^T\bfg]=\mathbb{E}_{P}\left[\sum_{k=1}^\np \sum_{j=1}^{\nh}\rowu_j^T g_j \mathbb{I}_{p_k=j} \right]=  \frac{\np}{\nh}\bfU^T \bfg = 0.
$$
}
\end{remark}
We now show that \gpod is robust with respect to noise in the sense that increasing the number of sampling points $\np$ reduces the effect of the noise.
\begin{thm}
Consider the same setup as in \Cref{lm:EigenBound} and \Cref{A:LEM:4.3}. Define
\begin{equation}  \label{definezeta}
\zeta =   \min\left\{  {\frac{1}{\sqrt{1 - \gamma}}}, \frac{1}{1-\gamma}\sqrt{\frac{\nd}{\np} \mu(\Ucal)}\,\right\}.
\end{equation}
Then,
\begin{equation}  \label{eqn:mainbound}
\mathbb{E}_{P}\left[\mathbb{E}_{\epsilon}\left[\|\bff - \shat{\bff}_{\bfepsilon}\|_2\right]\right] \leq
\left(1 + \zeta\right)\|\bff - \bfU\bfU^T\bff\|_2 + \frac{\|\bfsigma\|_\infty}{(1 - \gamma)} \sqrt{\frac{\nd\nh}{\np}}
\end{equation}
with probability at least $1 - \delta$, where the expectation $\mathbb{E}_P$ is with respect to the distribution of the samples and $\mathbb{E}_{\epsilon}$ with respect to the noise in $\bff_{\bfepsilon}$ as defined in \cref{eq:Prelim:NoisyF}.
\label{A:prop:MainBound}
\end{thm}
\begin{proof}
	We split the error following the strategy of \cite[Theorem~2]{Cohen2013}.
	With the triangular inequality, we obtain
	\begin{equation}
	\mathbb{E}_{P}\left[\mathbb{E}_{\epsilon}\left[\|\bff - \shat{\bff}_{\bfepsilon}\|_2\right]\right] \leq \mathbb{E}_P\left[\|\bff - \shat{\bff}\|_2\right] + \mathbb{E}_P\left[\mathbb{E}_{\bfepsilon}\left[\|\bfU(\bfP^T\bfU)^\dagger\bfP^T\bfepsilon\|_2\right]\right]\,.
	\label{A:eq:MainBound:FirstTriangle}
	\end{equation}
with $\shat{\bff} = \bfU(\bfP^T\bfU)^\dagger\bfP^T\bff$ and $\bff_{\bfepsilon} = \bff + \bfepsilon$.

	To bound the first term on the right-hand side of the inequality \cref{A:eq:MainBound:FirstTriangle}, set $\bfg = {\bff} - \bfU\bfU^T{\bff}$. Similarly to  the case in \cite{DEIM}, it holds that
\begin{equation}
\bff - \shat{\bff}= \bff - \bfU(\bfP^T\bfU)^\dagger\bfP^T\bff
 = \bfg - \bfU(\bfP^T\bfU)^\dagger\bfP^T\bfg\,,
\label{eq:fminusDf}
\end{equation}
where we used the fact  that  $(\bfP^T\bfU)^\dagger\bfP^T\bfU = \bfI_\nd $ with probability, at least, $1 - \delta$. Note that $\bfg = {\bff} - \bfU\bfU^T{\bff}$ is orthogonal to the range of $\bfU$. Then,
\Cref{A:LEM:4.3} implies that
\begin{equation}\label{eq:L8:bound}
	\mathbb{E}_{P}[\|\bfU(\bfP^T\bfU)^\dagger\bfP^T\bfg\|_2]
	\leq
  \min\left\{  {\frac{1}{\sqrt{1 - \gamma}}}, \frac{\sqrt{ \frac{\nd}{\np}\mu(\Ucal)}}{1-\gamma}\right\} \|\bfg\|_2
\end{equation}
with probability at least $1 - \delta$. Then, 
 \cref{eq:fminusDf} and the linearity of the expectation yield
\begin{equation}\label{eq:1+zeta}
\mathbb{E}_P\left[\|\bff - \shat{\bff}\|_2\right]  \leq \left(1 + \zeta\right)\|\bff - \bfU\bfU^T\bff\|_2,
\end{equation}
where $\zeta$ is as defined in \cref{definezeta}.
Now consider the second term on the right-hand side of \cref{A:eq:MainBound:FirstTriangle}.
Note that $\bfepsilon$ is not necessarily orthogonal to $\bfU$, and therefore \Cref{A:rem:gpU} cannot be applied. We make the approximations
\begin{equation}
\mathbb{E}_P\left[\mathbb{E}_{\bfepsilon}\left[\|\bfU(\bfP^T\bfU)^\dagger\bfP^T\bfepsilon\|_2\right]\right] \leq \frac{\nh}{(1 - \gamma)\np}\mathbb{E}_P\left[\mathbb{E}_{\bfepsilon}\left[\|(\bfP^T\bfU)^T\bfP^T\bfepsilon\|_2\right]\right]\,,
\label{A:eq:MainBound:MainApprox}
\end{equation}
which holds with probability at least $1 - \delta$, see \cref{ub2} and \cref{A:eq:DERM:LemmaPseudo:Inter2a} in the proof of \Cref{A:LEM:4.3}. Consider now $\mathbb{E}_P\left[\mathbb{E}_{\bfepsilon}\left[\|(\bfP^T\bfU)^T\bfP^T\bfepsilon\|_2^2\right]\right] $. With the same notation as in the proof of \Cref{A:LEM:4.3}, and building on the proof of \cite[Theorem~3]{Cohen2013}, we have
\begin{eqnarray*}
	&&	\mathbb{E}_P\left[\mathbb{E}_{\bfepsilon}\left[\|(\bfP^T\bfU)^T\bfP^T\bfepsilon\|_2^2\right]\right] = \mathbb{E}_P\left[\mathbb{E}_{\bfepsilon}\left[\bla (\bfP^T\bfU)^T\bfP^T\bfepsilon, (\bfP^T\bfU)^T\bfP^T\bfepsilon \bra_2 \right]\right] \\
	&& =\mathbb{E}_P\left[\mathbb{E}_{\bfepsilon}\left[\bla \sum_{i=1}^{\np} \epsilon_{p_i}\rowu_{p_i}^T, \sum_{j=1}^{\np} \epsilon_{p_j}\rowu_{p_j}^T  \bra_2 \right]\right] = \mathbb{E}_P\left[  \sum_{i=1}^m \sum_{j=1}^{\np} \mathbb{E}_{\bfepsilon}\left[  \bla \epsilon_{p_i}\rowu_{p_i}^T, \epsilon_{p_j} \rowu_{p_j}^T  \bra_2 \right]    \right] \\
	&& = \mathbb{E}_P\left[  \sum_{i=1}^{\np} \sum_{j=1}^{\np} \mathbb{E}_{\bfepsilon}\left[\epsilon_{p_i}\epsilon_{p_j}  \bla \rowu_{p_i}^T,  \rowu_{p_j}^T  \bra_2 \right]    \right]
	= \mathbb{E}_P\left[  \sum_{i=1}^{\np} \sum_{j=1}^{\np} \mathbb{E}_{\bfepsilon}\left[\epsilon_{p_i}\epsilon_{p_j}  \right] \bla \rowu_{p_i}^T,  \rowu_{p_j}^T  \bra_2   \right] \\
	&& = \mathbb{E}_P\left[  \sum_{i=1}^{\np}  \mathbb{E}_{\bfepsilon}\left[\epsilon_{p_i}^2  \right] \bla \rowu_{p_i}^T,  \rowu_{p_i}^T  \bra_2   \right] \;\;{\mbox{(since $\epsilon_{p_i}$ and $\epsilon_{p_j}$ are independent for $i\neq j$)}} \\
	&& = \mathbb{E}_P\left[  \sum_{i=1}^{\np}  \sigma_{p_i}^2  \bla \rowu_{p_i}^T,  \rowu_{p_i}^T  \bra_2   \right] =
	\mathbb{E}_P\left[  \sum_{i=1}^{\np} \sum_{j=1}^{\nh} \sigma_{j}^2   \| \rowu_{j}^T\|_2^2 \mathbb{I}_{p_i=j}  \right] = \frac{\np}{\nh}\sum_{j=1}^{\nh} \sigma_j^2 \| \rowu_{j}^T\|_2^2\,. \\
	\end{eqnarray*}
Using the fact that $\sum_{j=1}^{\nh}  \| \rowu_{j}^T\|_2^2 = \| \bfU\|_F^2 = n$, we obtain
$$
\mathbb{E}_P\left[\mathbb{E}_{\bfepsilon}\left[\|(\bfP^T\bfU)^T\bfP^T\bfepsilon\|_2^2\right]\right]  = \frac{\np}{\nh}\sum_{j=1}^{\nh} \sigma_j^2 \| \rowu_{j}^T\|_2^2\leq
\frac{\np}{\nh} \| \bfsigma\|_\infty^2 \sum_{j=1}^{\nh}  \| \rowu_{j}^T\|_2^2
= \frac{\np \nd }{\nh}\| \bfsigma\|_\infty^2.
$$
Applying Jensen's inequality, together with  \cref{A:eq:MainBound:MainApprox} yields
\begin{eqnarray}
\mathbb{E}_P\left[\mathbb{E}_{\bfepsilon}\left[\|\bfU(\bfP^T\bfU)^\dagger\bfP^T\bfepsilon\|_2\right]\right] &\leq &
\frac{\nh}{(1 - \gamma)\np} \nonumber
  \sqrt{\frac{\np\nd}{\nh}\|\bfsigma\|_\infty^2} = \frac{\|\bfsigma\|_\infty}{(1 - \gamma)}
  \sqrt{\frac{\nd\nh}{\np}}
\label{A:eq:MainBound:MainApproxC}
\end{eqnarray}
with probability at least $1 - \delta$. Combining this 
with \cref{eq:1+zeta} proves the theorem.
\end{proof}
\begin{remark}
\cool{\Cref{A:prop:MainBound} reveals that as $\np \to \infty$ (recall we perform uniform sampling with replacement), the upper bound in \cref{eqn:mainbound} converges to the
projection error
 $\| \bff - \bfU \bfU^T \bff\|_2$. In the numerical results in \Cref{sec:NumRes}, rather than investigating the asymptotic behavior as $\np \to \infty$, we will typically keep the ratio $\nd/\np$ low by increasing $\np$ with $\nd$ and so preventing that the noise term in \eqref{eqn:mainbound} dominates.}
 \end{remark}

 \section{The \odeime deterministic sampling strategy}
 \label{sec:DetSampling}
 In this section, we present a deterministic strategy that selects sampling points to reduce the quantity $\|(\bfP^T\bfU)^{\dagger}\|_2$, which controls how sensitive the \gpod oblique projection is to perturbations and noise, cf.~\Cref{sec:Prelim:Instability}. While our probabilistic analysis in \Cref{sec:ODEIM} shows that sampling points that are selected uniformly in $\{1, \dots, \nh\}$ are sufficient for \gpod to be robust with respect to noise, the number of uniformly selected sampling points that are required grows with, e.g., the coherence $\mu(\Ucal)$ of the space $\Ucal$. The following deterministic selection strategy aims to achieve robustness with fewer points than uniform sampling by taking information about the space $\Ucal$ into account. We refer to the Introduction in \Cref{sec:Introduction} for references to other sampling strategies.

 We propose the \odeime sampling algorithm that is based on lower bounds of the smallest eigenvalues of certain structured matrix updates introduced in \cite{IpsenBound} and that is a special case of the approach introduced by Zimmermann et al.~in \cite{ZimW16}. The ``E'' in \odeime stands for ``eigenvector''. The goal of the \odeime sampling algorithm is to select points that minimize $\|(\bfP^T\bfU)^{\dagger}\|_2$. This minimization problem is equivalent to maximizing the smallest singular value of $\bfP^T\bfU$ because
 \[
 \|(\bfP^T\bfU)^{\dagger}\|_2 = s_{\text{max}}\left((\bfP^T\bfU)^\dagger\right) = \frac{1}{s_{\text{min}}(\bfP^T\bfU)}\,,
 \]
 where $s_{\text{max}}(\bfM)$ and $s_{\text{min}}(\bfM)$ denote the largest and smallest singular value of the matrix $\bfM$, respectively. The \odeime algorithm relies on lower bounds of the smallest eigenvalues to select points that maximize $s_{\text{min}}(\bfP^T\bfU)$ by leveraging the eigenvector corresponding to the smallest eigenvalue.

 \subsection{Singular values after symmetric rank-one updates}
 Consider the basis matrix $\bfU$ and the sampling points matrix\footnote{Note that we have changed the notation slightly here  and added the subscript $``\np"$ to $\bfP$. This will help distinguish the sampling points matrix when new indices are added.} $\bfP_{\np}$ that takes $\np \geq \nd$ samples.
  Consider now the \textsf{SVD} of $\bfP_{\np}^T\bfU \in \mathbb{R}^{\np \times \nd}$
 \[
 \bfV_{\np}\bfSigma_{\np}\bfW_{\np}^T = \bfP_{\np}^T\bfU\,,
 \]
 where $\bfV_{\np} \in \mathbb{R}^{\np \times \nd}$ contains, as its columns, the left-singular vectors, the matrix
 $
 \bfSigma_{\np} = \operatorname{diag}[
 s_1^{(\np)},  \cdots,  s_{\nd}^{(\np)}
 ] \in \mathbb{R}^{\nd \times \nd}
 $
 is a diagonal matrix with the singular values $s_1^{(\np)}, \dots, s_{\nd}^{(\np)}$, in descending order, and $\bfW_{\np} \in \mathbb{R}^{\nd \times \nd}$ contains, as its columns, the right-singular vectors. Note that we assume that $\bfP_{\np}^T\bfU$ has full column rank in the following, which can be ensured by initializing \odeime with, e.g., the \qdeim interpolation points. If we add a sampling point, we obtain
 \[
 \bfP_{\np + 1}^T\bfU = \begin{bmatrix}
 \bfP_{\np}^T\bfU\\
 \bfu_{+}
 \end{bmatrix} \in \mathbb{R}^{\np + 1 \times \nd}\,,
 \]
 where $\bfu_+ \in \mathbb{R}^{1 \times \nd}$ is the row of $\bfU$ that is selected by the new sampling point. Following the work by Zimmermann et al.~\cite{ZimW16}, the change of the singular values of $\bfP_{\np}^T\bfU$ to $\bfP_{\np + 1}^T\bfU$ can be understood via a symmetric rank-one update.
 We have
 \[
 \bfP^T_{\np + 1}\bfU = \begin{bmatrix} \bfV_{\np} & 0 \\
 0 & 1
 \end{bmatrix} \begin{bmatrix}
 \bfSigma_{\np}\\
 \bfu_+ \bfW_{\np}
 \end{bmatrix}
 \bfW_{\np}^T\,.
 \]
 The singular values of $\bfP_{\np + 1}^T\bfU$ are given by the square roots of the eigenvalues of $(\bfP_{\np + 1}^T\bfU)^T(\bfP_{\np + 1}^T\bfU)$, which we represent as
 \[
 (\bfP_{\np + 1}^T\bfU)^T(\bfP_{\np + 1}^T\bfU) = \bfW_{\np} \left(\bfSigma_{\np}^2 + \bfW_{\np}^T\bfu_+^T\bfu_+\bfW_{\np}\right)\bfW_{\np}^T\,.
 \]
 Define $\bfLambda_{\np + 1} =\bfSigma_{\np}^2 + \bfW_{\np}^T\bfu_+^T\bfu_+\bfW_{\np}$.
 With $\bar{\bfu}_+ = \bfW_{\np}^T\bfu_+^T$, we obtain
 $\bfLambda_{\np + 1} = \bfSigma_{\np}^2 + \bar{\bfu}_+\bar{\bfu}_+^T$,
 which is a symmetric rank-one update to the diagonal matrix $\bfSigma_{\np}^2$. The square roots of the eigenvalues of $\bfLambda_{\np + 1}$ are the singular values of $\bfP^T_{\np + 1}\bfU$.

 Let $\lambda_1^{(\np)}, \dots, \lambda_{\nd}^{(\np)}$ be the eigenvalues of $\bfSigma_{\np}^2$ and let $\lambda_1^{(\np + 1)}, \dots, \lambda_{\nd}^{(\np + 1)}$ be the eigenvalues of $\Lambda_{\np + 1}$, both listed in descending order. Our goal is now to select a row of $\bfU$ that maximizes the smallest eigenvalue $\lambda_n^{(\np + 1)}$. From Weyl's theorem \cite{Weyl1912,WilkinsonBook} we have that $\lambda_{\nd}^{(\np + 1)} \geq \lambda_{\nd}^{(\np)}$, which shows that adding any sampling point will, at least, not increase $\|(\bfP_{\np + 1}^T\bfU)^\dagger\|_2$ compared to $\|(\bfP_{\np}^T\bfU)^\dagger\|_2$.

 \subsection{Lower bounds for eigenvalues of updated matrices}
 \label{sec:ODEIMDet:Bound}
 We now use the results by Ipsen et al.~in \cite{IpsenBound} to derive a heuristic strategy with the aim of selecting sampling points that lead to a fast increase of the smallest eigenvalue, i.e., to a fast decrease of $\|(\bfP_{\np}^T\bfU)^\dagger\|_2$.

 Let $g = \lambda_{\nd - 1}^{(\np)} - \lambda_{\nd}^{(\np)}$ be the eigengap. Note that we need $\lambda_{\nd - 1}^{(\np)} > \lambda_{\nd}^{(\np)}$ in the following. Let $\bfz_{\nd}^{(\np)} \in \mathbb{R}^{\nd}$ be the eigenvector of $\bfSigma_{\np}^2$ corresponding to the smallest eigenvalue $\lambda_{\nd}^{(\np)}$, with $\|\bfz_{\nd}^{(\np)}\|_2 = 1$. In our case $\bfz_{\nd}^{(\np)}$ is the $\nd$-th canonical unit vector of dimension $\nd$ because $\bfSigma_{\np}^2$ is diagonal with diagonal elements ordered descending. Then, as shown in \cite[Corollary~2.2]{IpsenBound},
 \begin{equation}
 \lambda_{\nd}^{(\np + 1)} \geq \lambda_{\nd}^{(\np)} + \frac{1}{2}\left(g + \|\bar{\bfu}_+\|_2^2 - \sqrt{\left(g + \|\bar{\bfu}_+\|_2^2\right)^2 - 4 g (\bfz_{\nd}^{(\np)}{}^T\bar{\bfu}_+)^2}\right)\,.
 \label{eq:IpsenBound}
 \end{equation}
 Observe that the bound \cref{eq:IpsenBound} depends on the eigen\emph{vector} corresponding to the smallest eigenvalue.

 The bound \cref{eq:IpsenBound} motivates us to add the rows of $\bfU$ that maximize
 \begin{equation}
 g + \|\bar{\bfu}_+\|_2^2 - \sqrt{\left(g + \|\bar{\bfu}_+\|_2^2\right)^2 - 4 g (\bfz_{\nd}^{(\np)}{}^T\bar{\bfu}_+)^2}\,.
 \label{eq:ODEIM:OurCriterion}
 \end{equation}
 The criterion \cref{eq:ODEIM:OurCriterion} is related to the criteria developed in \cite{ZimW16}. While we build on the perturbation bounds introduced in \cite{IpsenBound}, the authors of \cite{ZimW16} directly derive criteria that take the eigenvector $\bfz_{\nd}$ corresponding to the smallest eigenvalue into account, see \cite[page~A2834]{ZimW16} and \cite[Remark~2, item~3]{ZimW16}. In fact, the work \cite{ZimW16} goes a step further and also takes into account inner products with eigenvectors corresponding to larger eigenvalues. We do not consider these additional steps discussed in \cite{ZimW16} in the following.

 \subsection{The \odeime algorithm}
 The \odeime sampling approach that we consider is summarized in \Cref{alg:ODEIMOversampling}. It iteratively selects new sampling points that maximize \cref{eq:ODEIM:OurCriterion} in a greedy fashion. In line~2 of \Cref{alg:ODEIMOversampling}, the first $\nd$ points are selected with \qdeim, see \Cref{alg:qdeim}. Then, for each point $i = \nd + 1, \dots, \np$, the \textsf{SVD} of $\bfP_i^T\bfU$ (which is \texttt{U(p, :)} in the notation used in \Cref{alg:ODEIMOversampling}) is computed to obtain the right-singular vectors as columns of the matrix $\bfW_{m}$. The eigengap $g$ is computed on line 6 in \Cref{alg:ODEIMOversampling}. Then, $\bar{\bfU} = \bfW_i^T\bfU^T$ is obtained on line 7 in \Cref{alg:ODEIMOversampling}. The bound \eqref{eq:ODEIM:OurCriterion} is then computed from $\bar{\bfU}$ for each column $\bar{\bfu}_+$ in lines 8--9 and sorted descending on line 10. On line 9, it is exploited that the eigenvector $\bfz_{\nd}^{(i)}$ is the $\nd$-th canonical unit vector of dimension $\nd$ and so $\bfz_{\nd}^{(i)}{}^T\bar{\bfu}_+$ in \eqref{eq:ODEIM:OurCriterion} for column $\bar{\bfu}_+$ of $\bar{\bfU}$ is computed as $\bfz_{\nd}^{(i)}{}^T\bar{\bfu}_+ = \bfz_{\nd}^{(i)}\bfW_i^{T}\bfu_+^T = (\bfw_i^{\text{end}})^T\bfu_+^T$, where $\bfw_i^{\text{end}}$ is the right-singular vector corresponding to the smallest singular value. This means that $(\bfz_{\nd}^{(i)})^T\bar{\bfU}$  is given by the last row of $\bar{\bfU}$ (denoted as \texttt{Ub(end, :)} in \Cref{alg:ODEIMOversampling}). The point corresponding to the column of $\bar{\bfU}$ (row of $\bfU$) with the largest value \eqref{eq:ODEIM:OurCriterion} is added as sampling point and the procedure is repeated.
 Each iteration in \odeime requires performing an \textsf{SVD} of a small matrix whose size grows with the reduced dimension $\nd$ and the number of sampling points. Each \textsf{SVD} is in $\mathcal{O}(\nd^2\np)$ (for $\np > \nd$). Thus, selecting $\np$ points with \odeime is in $\mathcal{O}(\nd^2\np^2)$. Note that the sampling point selection is performed during the construction of the reduced model in the offline phase.

 \cool{
 The \odeime algorithm returns points that are not necessarily nested with respect to the dimension of the \deim basis: Consider a basis matrix $\bfU_n$ with $n$ columns and the corresponding set $\mathcal{P}_n$ of $m > n$ points selected by \odeime. Let now $\bfU_{n + 1}$ be a basis matrix with $n + 1$ columns where the first $n$ columns coincide with the columns of $\bfU_n$ and let $\mathcal{P}_{n + 1}$ be the set of at least $m$ points selected with \odeime. Then, it is possible that $\mathcal{P}_n \not\subset \mathcal{P}_{n + 1}$, which is in contrast to, e.g., the greedy EIM algorithm \cite{barrault04-EIM}, for which $\mathcal{P}_n \subset \mathcal{P}_{n + 1}$ holds if $n$ and $n + 1$ points are selected, respectively. Nestedness of points is a desired property in situations where one wants to, for example, rapidly and adaptively select the number of \deim basis vectors and sampling points without running the sampling algorithm from scratch. One such situation is in model reduction when the dimension of the \deim space is selected during the online phase. One option to avoid running \odeime during the online phase in this situation is to pre-compute \odeime sampling points for a range of dimensions of the \deim space and to store them during the offline phase. In the online phase, the pre-computed points can be quickly loaded depending on the basis dimension that is selected online, instead of running \odeime during the online phase.
 }

 \begin{algorithm}[t]
 \caption{Sampling points selection with \odeime (Matlab notation)}\label{alg:ODEIMOversampling}
 \begin{lstlisting}[style=Matlab-editor,mathescape]
 function [ p ] = gpode( U, m )
 [~, ~, p] = qr(U', 'vector');
 p = p(1:size(U, 2))';
 for i=length(p)+1:m
     [~, S, W] = svd(U(p, :), 0);
     g = S(end-1,$~$end-1).^2 - S(end,$~$end)^2;
     Ub = W'*U';
     r = g + sum(Ub.^2, 1);
     r = r-sqrt((g+sum(Ub.^2,1)).^2-4*g*Ub(end, :).^2);
     [~, I] = sort(r, 'descend');
     e = 1;
     while any(I(e) == p)
         e = e + 1;
     end
     p(end + 1) = I(e);
 end
 end
 \end{lstlisting}

 \end{algorithm}

 \section{Numerical results}
 \label{sec:NumRes}
 This section compares the stability of \deim and \gpod with randomized and deterministic sampling algorithms on numerical examples. \Cref{sec:NumRes:Toy} revisits the toy example from \Cref{sec:Prelim:Instability} and demonstrates that \gpod provides stable approximations compared to \deim. \Cref{sec:NumRes:Combustion} approximates velocity fields from noisy measurements of single-injector combustion processes following the procedure introduced in \cite{WILLCOX2006208,BuiDW}. \Cref{sec:NumRes:Laplace} demonstrates the effect of taking more sampling points than basis vectors on a diffusion-reaction problem, where \gpod provides stable approximations in contrast to \deim.

 \subsection{Synthetic example}
 \label{sec:NumRes:Toy}
\begin{figure}
 \centering
 \begin{tabular}{cc}
 \multicolumn{2}{c}{{\LARGE\resizebox{0.95\columnwidth}{!}{\input{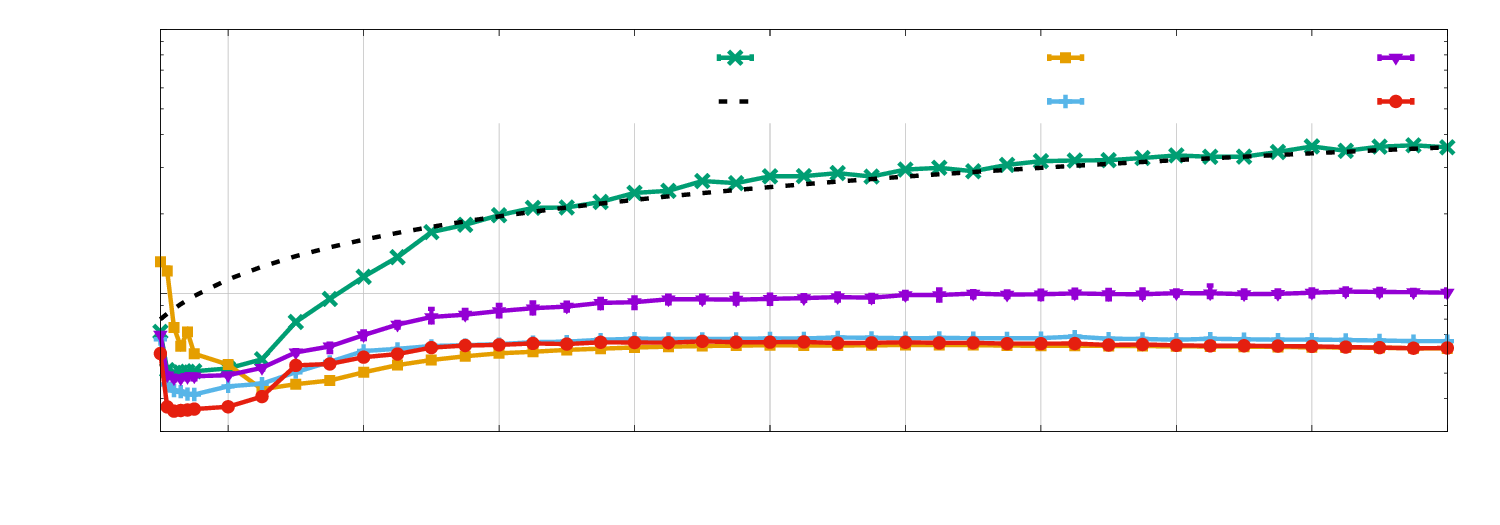}}}}\\
 \multicolumn{2}{c}{\scriptsize (a) std.~deviation $\sigma = 10^{-4}$}\\
 {\LARGE\resizebox{0.46\columnwidth}{!}{\input{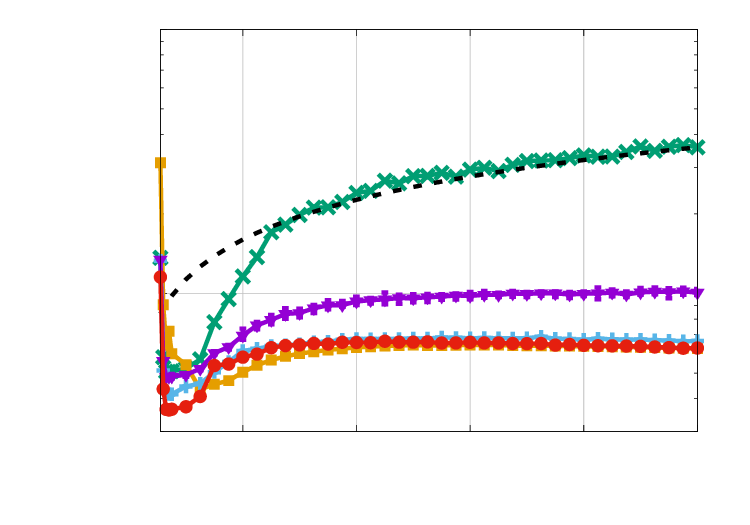}}} &  {\LARGE\resizebox{0.46\columnwidth}{!}{\input{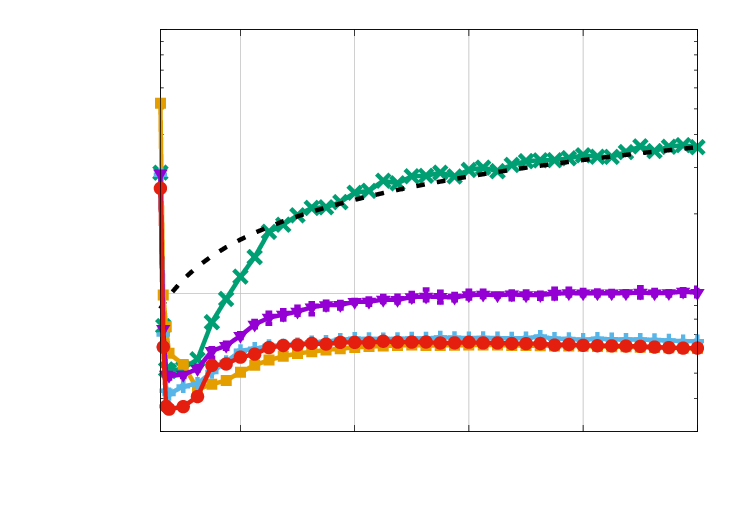}}}\\
 \scriptsize (b) std.~deviation $\sigma = 10^{-5}$ & \scriptsize (c) std.~deviation $\sigma = 10^{-6}$\\
 {\LARGE\resizebox{0.46\columnwidth}{!}{\input{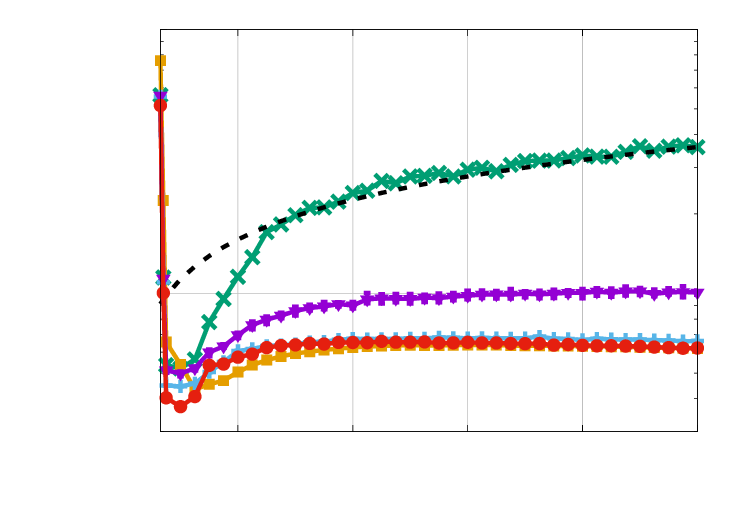}}} &  {\LARGE\resizebox{0.46\columnwidth}{!}{\input{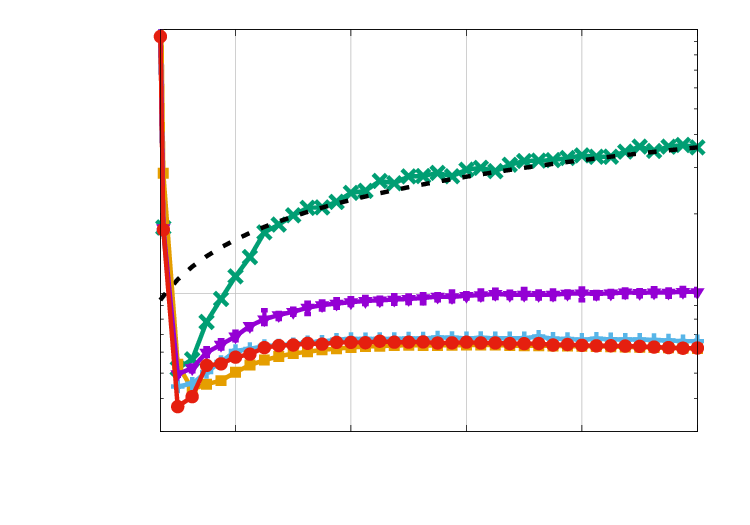}}}\\
 \scriptsize (d) std.~deviation $\sigma = 10^{-7}$ & \scriptsize (e) std.~deviation $\sigma = 10^{-8}$
 \end{tabular}
 \caption{Synthetic example: \gpod with more sampling points than basis vectors shows stable behavior and avoids the amplification of the error with a rate $\sqrt{\nd}$ as observed in the \qdeim approximation. All sampling strategies for \gpod give stable approximations in this example.}
 \label{fig:DEIMToyOver}
 \end{figure}

 Let us revisit the synthetic example introduced in \Cref{sec:Prelim:Instability}. We use the same setup as before but now approximate the noisy function with \gpod that takes more sampling points than basis vectors. We compare our \odeime sampling strategy to three sampling strategies from the literature. First, there is \odeimr, which takes the first $\nd$ sampling points with \qdeim and the subsequent $\np - \nd$ sampling points uniform randomly with replacement. Thus, the ``\textsf{R}'' in \odeimr stands for ``random''. Second, the strategy \odeiml takes the first $\nd$ points with \qdeim and the subsequent $\np - \nd$ points based on leverages scores as described in, e.g., \cite[Section~V.B]{AstWWB04}. Thus, the ``\textsf{L}'' in \odeiml stands for ``leverage scores''. Third, with \odeimd we denote the sampling strategy introduced in \cite[Algorithm~4]{Carlberg2013} that selects $\np > \nd$ sampling points by extending the \deim greedy algorithm \cite{Barrault2004,DEIM}. Thus, the ``\textsf{D}'' in \odeimd stands for ``\deim greedy''. The number of sampling points is set to $\np = 2\nd$ in case of \gpod in the following. We perform 10 replicates of the experiments and compute  the error as defined in \cref{eq:Prelim:L2ErrorNoise}. The results in Figure~\ref{fig:DEIMToyOver}a indicate that \gpod with more sampling points than basis vectors avoids the unstable behavior obtained with \deim, as suggested by our analysis presented in Section~\ref{sec:ODEIM}. All sampling algorithms perform well in this example, with \odeime, \odeiml, \odeimd achieving the lowest errors. Similar results are obtained for $\sigma \in \{10^{-5}, \dots, 10^{-8}\}$ in this example as shown in Figure~\ref{fig:DEIMToyOver}b-e. Note that the error decays linearly with $\sigma$ as long as the noise limits the overall accuracy rather than the projection error, as indicated by Theorem~\ref{A:prop:MainBound}. The error bars in Figure~\ref{fig:DEIMToyOver} show the minimum and maximum of the error over the 10 replicates. Except for \odeimr in Figure~\ref{fig:DEIMToyOver}a near $\nd = 1000$, the error bars are barely visible in the plots, which indicates that even small perturbations due to noise lead to unstable behavior in \qdeim.

\begin{figure}
 \centering
 \begin{tabular}{cc}
 {\LARGE\resizebox{0.46\columnwidth}{!}{\input{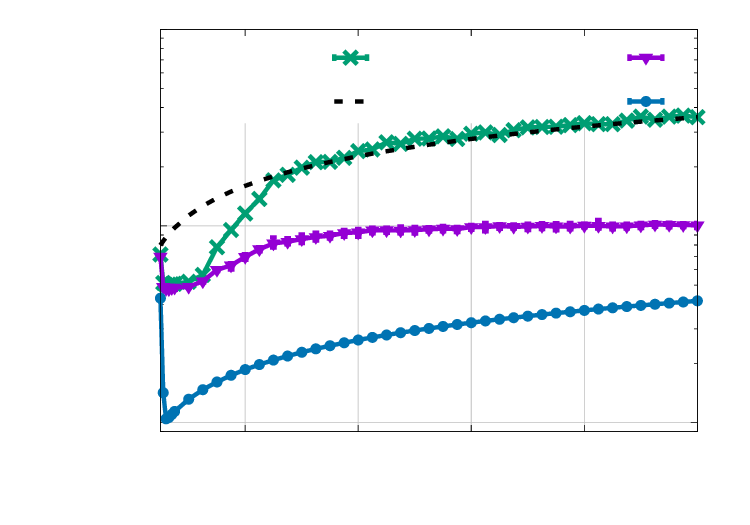}}} &  {\LARGE\resizebox{0.46\columnwidth}{!}{\input{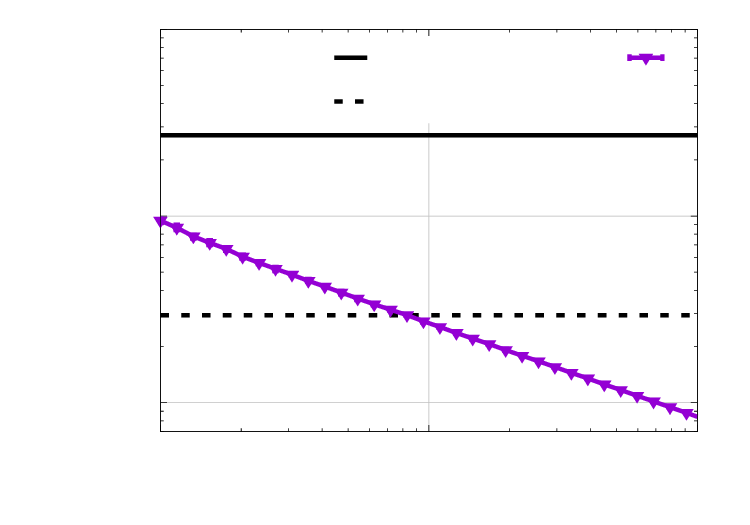}}}\\
 \scriptsize (a) sampling with ratio $\np/\nd = 2$ & \scriptsize (b) dimension $\nd = 500$ fixed
 \end{tabular}
 \caption{\cool{Plot (a) shows the projection error \eqref{eP} of noisy data with std.~deviation $\sigma = 10^{-4}$. The projection error \eqref{eP} grows with the space dimension $\nd$ because the Gaussian random noise vector can be more accurately approximated in the space as the dimension is increased. Plot (b) indicates that the error \eqref{eq:Prelim:L2ErrorNoise} of \odeimr converges with a rate $\sqrt{1/\np}$ to the projection error of noise-less data \eqref{eq:NumRes:ProjError} for a fixed dimension, cf.~Theorem~\ref{A:prop:MainBound}. The projection error of noise-less data is below machine precision in this example.}}
 \label{fig:DEIMToyOver2}
 \end{figure}

 \cool{Figure~\ref{fig:DEIMToyOver2}a compares the orthogonal projection error of the noisy data,
 \begin{equation}
 \sum_{j = 1}^k \frac{1}{2500}\sum_{i = 1}^{2500} \frac{\|\bff(\bfx; \xi_i^{\prime}) - \bfU\bfU^T\bff_{\bfepsilon_j}(\bfx; \xi_i^{\prime})\|_2}{\|\bff(\bfx; \xi_i^{\prime})\|_2},
 \label{eP}
 \end{equation}
 to error \eqref{eq:Prelim:L2ErrorNoise} of \qdeim and \odeimr. Note that the projection error \eqref{eP} grows with the dimension $\nd$ since the Gaussian random noise vector can be better and better approximated in the subspace spanned by the columns of $\bfU$ as the dimension $\nd$ is increased.  Figure~\ref{fig:DEIMToyOver2}b shows that if the dimension $\nd = 500$ is fixed and the number of sampling points $\np$ is increased, then the error \eqref{eq:Prelim:L2ErrorNoise} of \odeimr decays with a rate $\sqrt{1/\np}$ to the projection error of noise-less data, i.e,,
 \begin{equation}
 \sum_{j = 1}^k \frac{1}{2500}\sum_{i = 1}^{2500} \frac{\|\bff(\bfx; \xi_i^{\prime}) - \bfU\bfU^T\bff(\bfx; \xi_i^{\prime})\|_2}{\|\bff(\bfx; \xi_i^{\prime})\|_2}\,,
 \label{eq:NumRes:ProjError}
 \end{equation}
 which is below machine precision in this example for $\nd = 500$. The results shown in Figure~\ref{fig:DEIMToyOver2}b are in alignment with Theorem~\ref{A:prop:MainBound}, which states that the error of \gpod with uniform sampling with replacement converges with rate $\sqrt{1/\np}$ to the projection error.}

 \subsection{Velocity field approximations from noisy measurements of single-injector combustion process}
 \label{sec:NumRes:Combustion}
\begin{figure}
 \centering
 \resizebox{0.9\columnwidth}{!}{\input{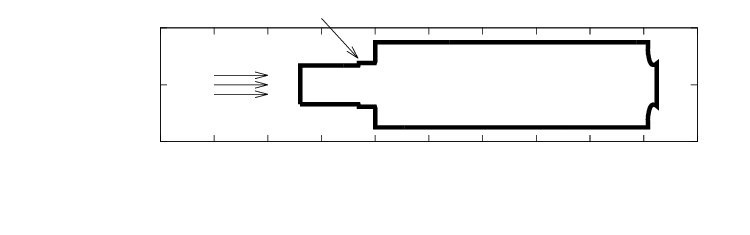}}
 \caption{Combustion: Geometry of combustion chamber; see, \cite{doi:10.1063/1.4916673,doi:10.2514/6.2019-2012} for details}
 \label{fig:NumRes:CombGeo}
 \end{figure}

 We consider velocity field approximations from noisy measurements of the single-injector combustion process described in detail in \cite{swischuk2019learning}. The combustion model follows the implementation of the General Equation and Mesh Solver (GEMS) code \cite{doi:10.1063/1.4916673,doi:10.2514/6.2019-2012} developed by Purdue University. The domain of the setup of \cite{swischuk2019learning} is shown in Figure~\ref{fig:NumRes:CombGeo}. Fuel and oxidizer are input with constant mass flow rates of $5.0 \frac{\text{kg}}{s}$ and $0.37 \frac{\text{kg}}{s}$, respectively. The fuel is composed of gaseous methane and the oxidizer is 42\% gaseous $\text{O}_2$ and 58\% gaseous $\text{H}_2\text{O}$. Details of the physics of the problem setup are described in \cite{doi:10.2514/6.2019-2012}.

\begin{figure}
 \begin{tabular}{cc}
 \multicolumn{2}{c}{\LARGE\resizebox{0.9\columnwidth}{!}{\input{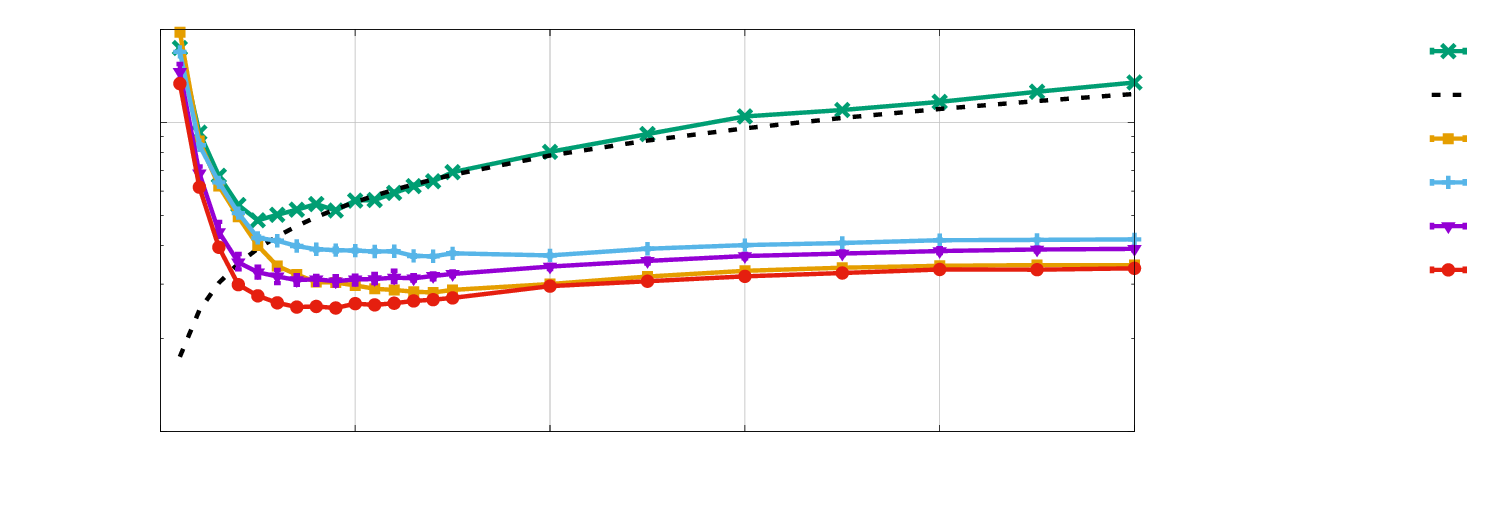}}}\\
 \multicolumn{2}{c}{\scriptsize (a) std.~deviation $\sigma = 3.4$ (about $1\%$ noise)}\\
 \LARGE\resizebox{0.45\columnwidth}{!}{\input{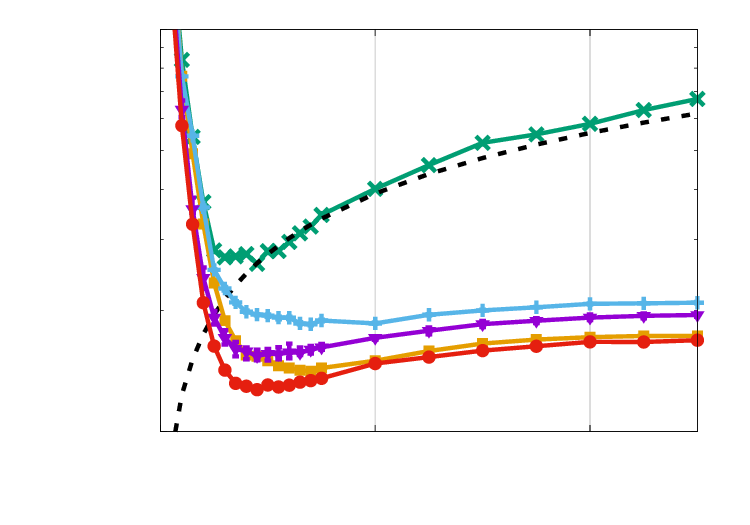}} &\LARGE\resizebox{0.45\columnwidth}{!}{\input{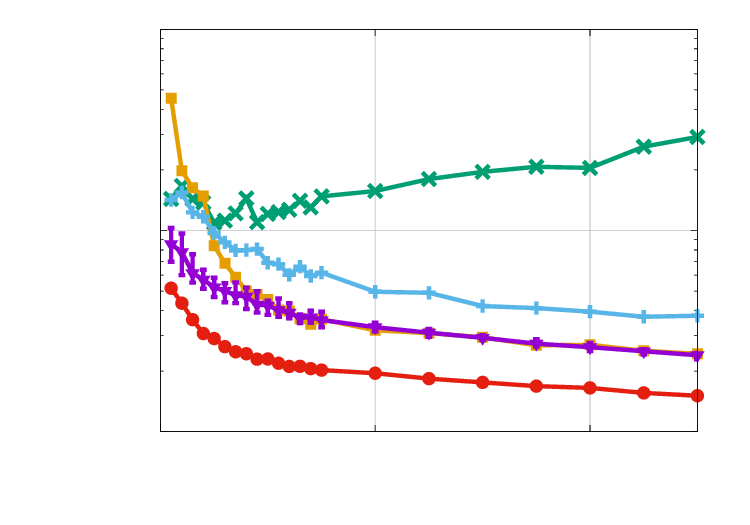}}\\
 \scriptsize (b) std.~deviation $\sigma = 1.7$ (about $0.5\%$ noise) & \scriptsize (c) norm of sampling operator\\
 \end{tabular}
 \caption{Combustion: The plots in this figure show that approximating the velocity field of the single-injector combustion process considered in this example from noisy measurements suffers from the instability described in Section~\ref{sec:Prelim:Instability} if \qdeim is used with the same number of sampling points as the dimension of the reduced space. In contrast, \gpod with various sampling strategies yields stable approximations, i.e., avoids the growth with rate $\sqrt{\nd}$ with the dimension $\nd$ of the reduced space. \odeime achieves the lowest error in this example.}
 \label{fig:Combustion}
 \end{figure}

 To generate snapshot data, the GEMS code is used to simulate the system for 0.7ms with a time step size of $\delta t = 10^{-7}$. The simulation leads to 7000 snapshots $\bfx(t_1), \dots, \bfx(t_{7000}) \in \mathbb{R}^{77046}$ at 7000 time points $t_1, \dots, t_{7000}$. The snapshots are of length $77046$ (there are 38523 spatial discretization points) and contain the velocity in x and y directions. The basis matrix $\bfU$ is derived with POD from the snapshots, where every fourth snapshot is skipped and kept as a test snapshot. Thus, the basis matrix is constructed from 5250 snapshots that are the columns of $\bfX$ and the 1750 test snapshots are ignored during construction of the basis matrix and collected as columns in $\bfX^{(\text{test})}$. The test snapshots are polluted with zero-mean Gaussian noise with standard deviation $\sigma = 1.7$ and $\sigma = 3.4$, respectively, and collected in the noisy test snapshot matrix $\bfX^{(\text{test})}_{\bfepsilon}$. A standard deviation of $\sigma = 1.7$ corresponds to about 0.5\% noise with respect to the mean of the snapshot matrix. Correspondingly, $\sigma = 3.4$ means that about 1\% noise is added to the test snapshots. Figure~\ref{fig:Combustion}a shows for $\sigma = 3.4$ the relative state error
 \begin{equation}
 \sum_{j = 1}^{10}\frac{\|\bfX^{(\text{test})} - \bfU(\bfP^T\bfU)^{\dag}\bfP^T\bfX^{(\text{test})}_{\bfepsilon_j}\|_F}{\|\bfX^{(\text{test})}\|_F}\,,
 \label{eq:Combustion:Error}
 \end{equation}
 over $j = 1, \dots, 10$ replicates of noise $\bfepsilon_j$. The matrix $\bfP$ is derived from \qdeim with $\np = \nd$ and from \odeimd, \odeiml, \odeimr, \odeime, respectively, with $\np = 2\nd$, i.e., twice as many sampling points as number of basis vectors. Figure~\ref{fig:Combustion}a shows the growth of the error \eqref{eq:Combustion:Error} for \qdeim, which uses the same number of sampling points as basis vectors. In contrast, taking more sampling points than basis vectors with $\gpod$ yields a stable approximation. All sampling strategies help to reduce the error \eqref{eq:Combustion:Error} significantly, where $\odeime$ achieves the lowest error in this example. Error bars show the minimum and the maximum of the error over the replicates. The error bars are barely visible in Figure~\ref{fig:Combustion}a, which indicates that small perturbations can lead to unstable behavior in \qdeim and that \gpod with more sampling points than basis vectors robustly gives stable approximations. Similar results are obtained for $\sigma = 1.7$ in Figure~\ref{fig:Combustion}b. Figure~\ref{fig:Combustion}c shows that the norm of the sampling operator $\|(\bfP^T\bfU)^{\dag}\|_2$ is lowest for $\odeime$, which is in agreement with the results in Figure~\ref{fig:Combustion}a and Figure~\ref{fig:Combustion}b.

 \subsection{Diffusion-reaction problem with nonlinear reaction term}
 \label{sec:NumRes:Laplace}
 We now demonstrate the stability of \gpod and \deim on a reduced model of a diffusion-reaction problem. The example demonstrates that instabilities in the \deim approximations can lead to unstable reduced models, which can be avoided with \gpod if more sampling points than basis vectors are used.

 \subsubsection{Problem setup}
 Let $\Omega = (0, 1)^2 \subset \mathbb{R}^2$ and consider the PDE
 \begin{equation}
 -\Delta u(\bfomega; \bfxi) + f(u(\bfomega; \bfxi); \bfxi) = 100\sin(2\pi \omega_1)\sin(2\pi \omega_2)\,,\qquad \bfomega \in \Omega\,,
 \label{eq:NLExp:PDE}
 \end{equation}
 where $\bfomega = [\omega_1, \omega_2]^T$ is the spatial coordinate, $u: \Omega \times \Dcal \to \mathbb{R}$ is the solution function, and $f: \mathbb{R} \times \Dcal \to \mathbb{R}$ is a nonlinear function
 \[
 f(u; \bfxi) = (0.1\sin(\xi_1)+2)\exp(-2.7\xi_1^2)(\exp(\xi_2 u 1.8) - 1)\,,
 \]
 with parameter $\bfxi = [\xi_1, \xi_2]^T \in \Dcal$. The PDE \cref{eq:NLExp:PDE} is closed with homogeneous Dirichlet boundary conditions. This example is a modification of the example considered in \cite{grepl_maday_nguyen_patera_2007}.

 \cool{We discretize \cref{eq:NLExp:PDE} with a second-order finite difference scheme on an equidistant mesh with mesh width $h = 1/255$ in $\Omega$, which leads to the state dimension $\nh = 65536$.} The system of nonlinear equations is solved with Newton's method. 
 We derive a reduced model from 1600 snapshots corresponding to a $40 \times 40$ grid of parameter values in the domain $\Dcal = [-\pi/2, \pi/2] \times [1, 5]$. The grid is equidistant in the first direction and logarithmically equidistant in the second direction. The basis matrix $\bfV$ is constructed with \rpod. The \rpod dimension is chosen as $\nr = 50$. The nonlinear term is approximated with empirical interpolation, with more details to follow below.
 The reduced model is tested on parameters corresponding to the $9 \times 9$ grid in $\Dcal$ that is linearly equidistant in the first direction and logarithmically equidistant in the second direction. The full-model states corresponding to the test parameters are collected in the test snapshot matrix $\bfX^{(\text{test})} \in \mathbb{R}^{\nh \times 81}$.

 \subsubsection{Results}
 We compare reduced models that differ in the way the nonlinear term is approximated. With ``\qdeim'' we denote the reduced models that approximate the nonlinear terms with \qdeim, which takes $\np = \nd$ sampling points. Reduced models that approximate the nonlinear term with \gpod are denoted as ``\odeimd'', ``\odeimr'', ``\odeiml'', and ``\odeime'', respectively, depending on which sampling strategy is used.

 \begin{figure}
 \centering
 \begin{tabular}{cc}
 {\Large\resizebox{0.46\columnwidth}{!}{\input{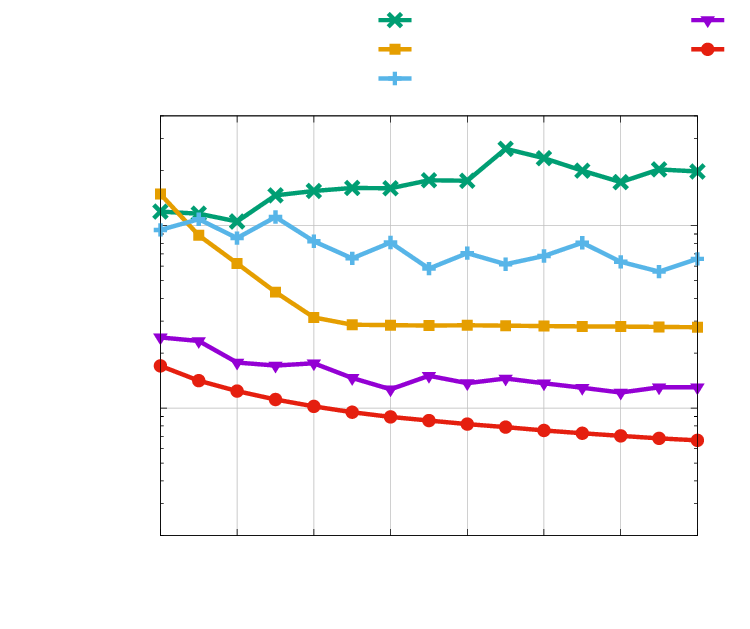}}} & {\Large\resizebox{0.46\columnwidth}{!}{\input{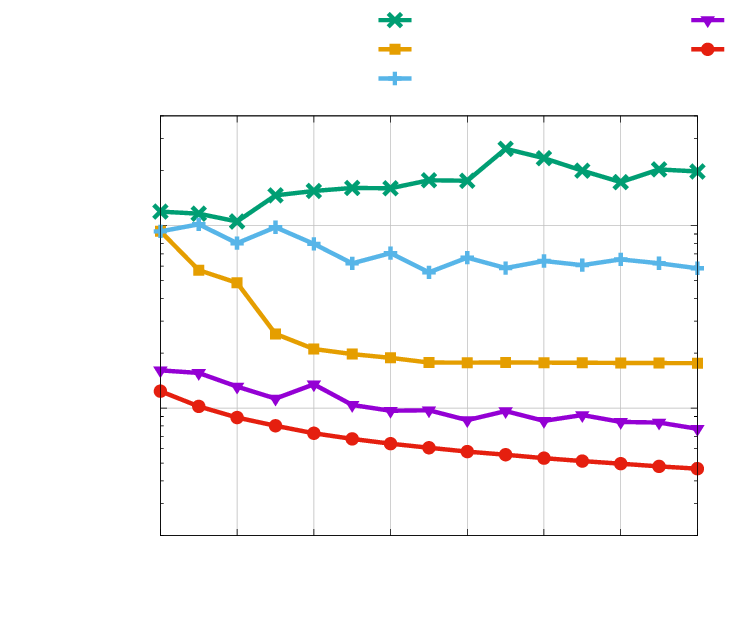}}}\\
 \scriptsize (a) $\np = 4\nd$ (factor 4 oversampling) & \scriptsize (b) $\np = 8\nd$ (factor 8 oversampling)\\
 \end{tabular}
 \caption{Diffusion reaction example: The sampling operators derived with \odeime achieves the lowest norm $\|(\bfP^T\bfU)^{\dag}\|_2$ in this example.}
 \label{fig:NLExpLConst}
 \end{figure}

\begin{figure}
 \centering
 \begin{tabular}{cc}
 {\LARGE\resizebox{0.47\columnwidth}{!}{\input{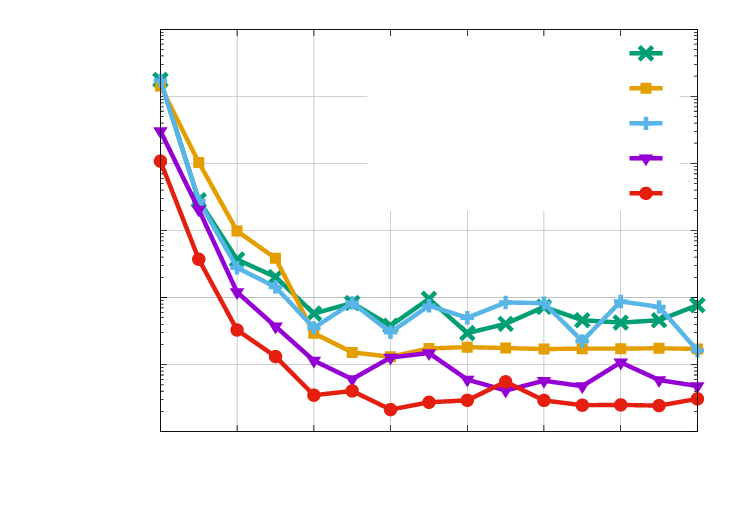}}} & \hspace*{-0.5cm}{\LARGE\resizebox{0.47\columnwidth}{!}{\input{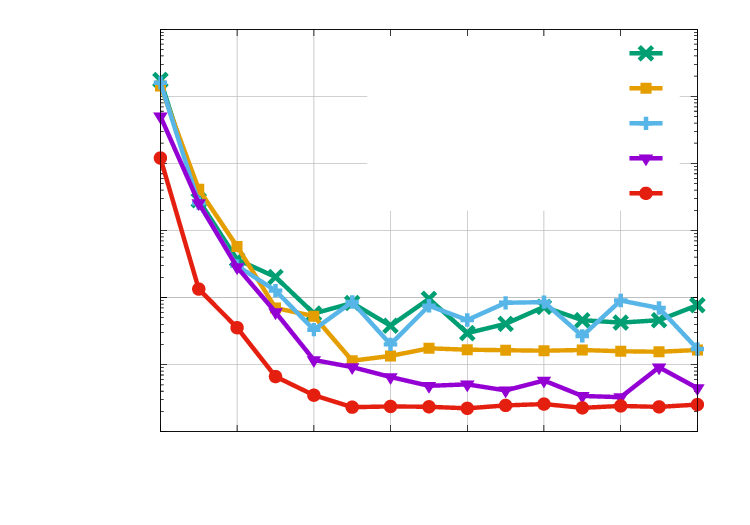}}}\\
 \scriptsize (a) oversampling $\np/\nd = 4$ & \scriptsize (b) oversampling $\np/\nd = 8$\\
 \end{tabular}
 \caption{Diffusion reaction example: Without noise, \qdeim and \gpod show stable behavior in this example. Note that the dimension $\nr$ of the POD space is fixed and therefore the curves level off even if the dimension $\nd$ of the reduced space spanned by the columns of $\bfU$ is increased.}
 \label{fig:NLExpAllNoNoise}
 \end{figure}

 Figure~\ref{fig:NLExpLConst} compares the norm of the sampling operators for $\np = 4\nd$ and $\np = 8\nd$ for dimensions $\nd \in \{50, \dots, 400\}$. \odeime  \cool{ provides the sampling operator} with the lowest norm  in this example. We first run the reduced models for the test parameters without adding noise and collect the corresponding states as columns in the matrix $\tilde{\bfX}^{(\textsf{S})} \in \mathbb{R}^{\nr \times 81}$ where \textsf{S} is either \qdeim, \odeimd, \odeiml, \odeimr, or \odeime. The averaged relative state error
 \begin{equation}
 \frac{\|\bfX - \bfV \tilde{\bfX}^{(\textsf{S})}\|_F}{\|\bfX\|_F}
 \label{eq:NLExp:ErrorNoiseless}
 \end{equation}
 is shown in Figure~\ref{fig:NLExpAllNoNoise}a for an oversampling factor $\np/\nd = 4$ and in Figure~\ref{fig:NLExpAllNoNoise}b for $\np/\nd = 8$ for \textsf{S} either \qdeim, \odeimd, \odeiml, \odeimr, or \odeime. \qdeim as well as \gpod with all sampling strategies achieve stable approximations in the sense described in Section~\ref{sec:Prelim:Instability}, i.e., the error \eqref{eq:NLExp:ErrorNoiseless} does not grow with the dimension $\nd$ of the reduced space.

 \begin{figure}
 \centering
 {\LARGE\resizebox{0.80\columnwidth}{!}{\input{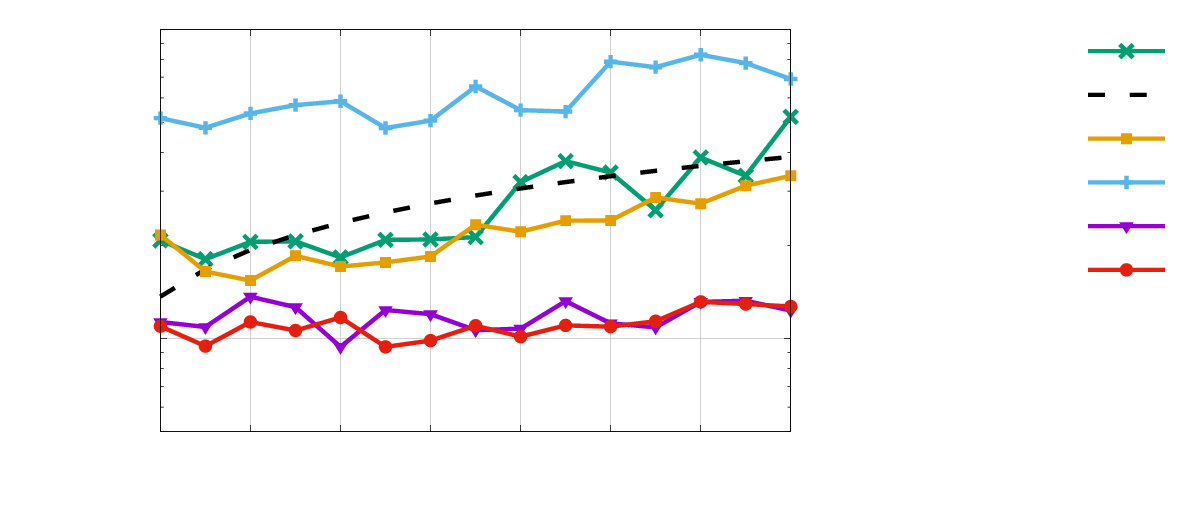}}}\\
 \caption{Diffusion reaction example \cool{with noise}: Approximating the nonlinear terms in this example with \qdeim leads to unstable behavior, which is indicated in this plot with a growth with rate $\sqrt{\nd}$. \gpod with more sampling points than basis vectors leads to stable reduced models with sampling strategies \odeime and \odeimr in this example.  Standard deviation of noise is $\sigma = 10^{-2}$ and oversampling factor is $\np/\nd = 4$.}
 \label{fig:NLExpAll}
 \end{figure}

\begin{figure}
 \centering
 \begin{tabular}{cc}
 {\Large\resizebox{0.46\columnwidth}{!}{\input{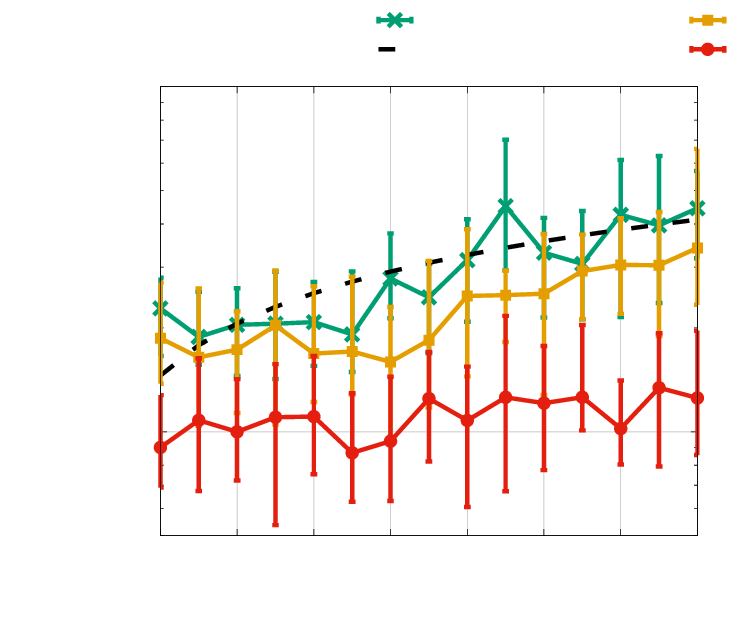}}} & {\Large\resizebox{0.46\columnwidth}{!}{\input{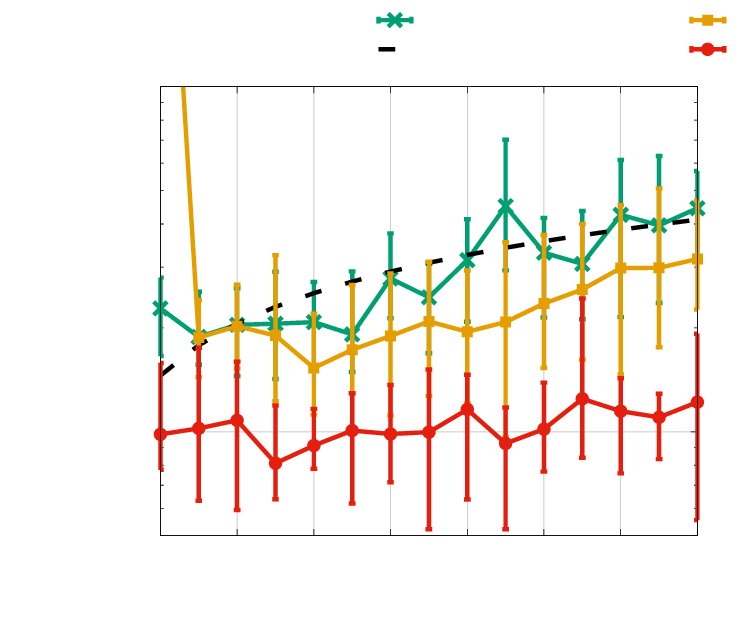}}}\\
 \scriptsize (a) std.~deviation $\sigma = 10^{-2}$, oversampling $\np = 4\nd$ & \scriptsize \hspace*{-0.5cm}(b) std.~deviation $\sigma = 10^{-2}$, oversampling $\np = 8\nd$\\
 {\Large\resizebox{0.46\columnwidth}{!}{\input{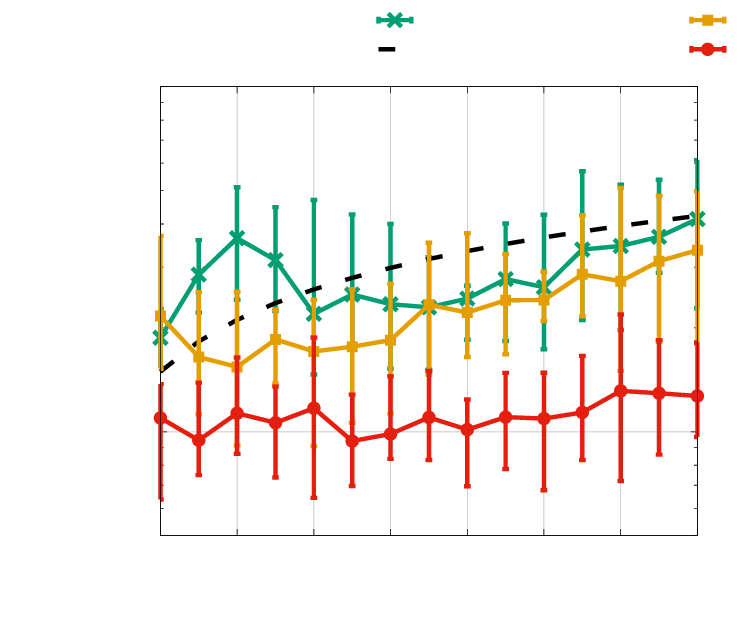}}} & {\Large\resizebox{0.46\columnwidth}{!}{\input{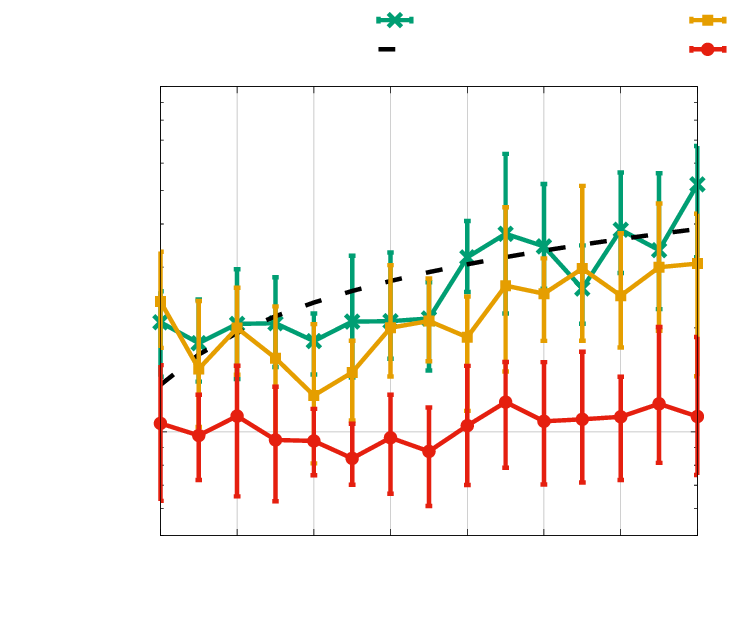}}}\\
 \scriptsize (c) std.~deviation $\sigma = 10^{-3}$, oversampling $\np = 4\nd$ & \scriptsize \hspace*{-0.5cm}(d) std.~deviation $\sigma = 10^{-3}$, oversampling $\np = 8\nd$\\
 {\Large\resizebox{0.46\columnwidth}{!}{\input{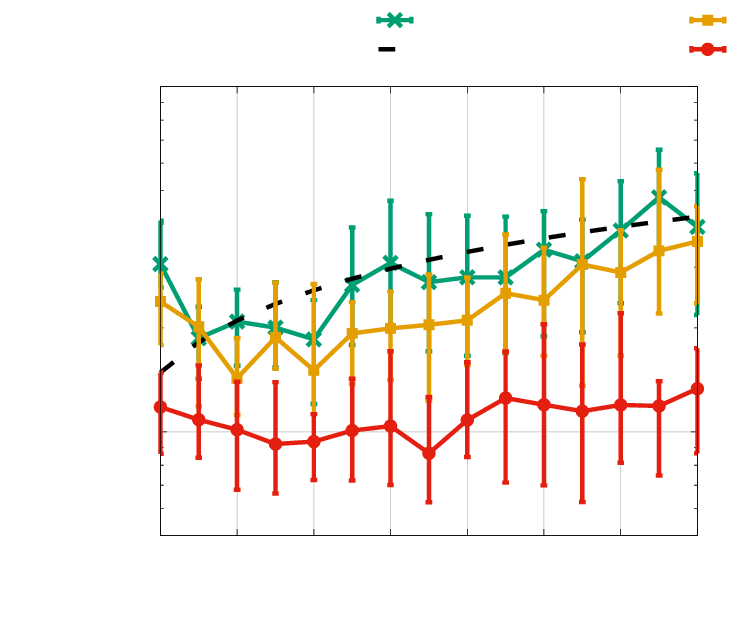}}} & {\Large\resizebox{0.46\columnwidth}{!}{\input{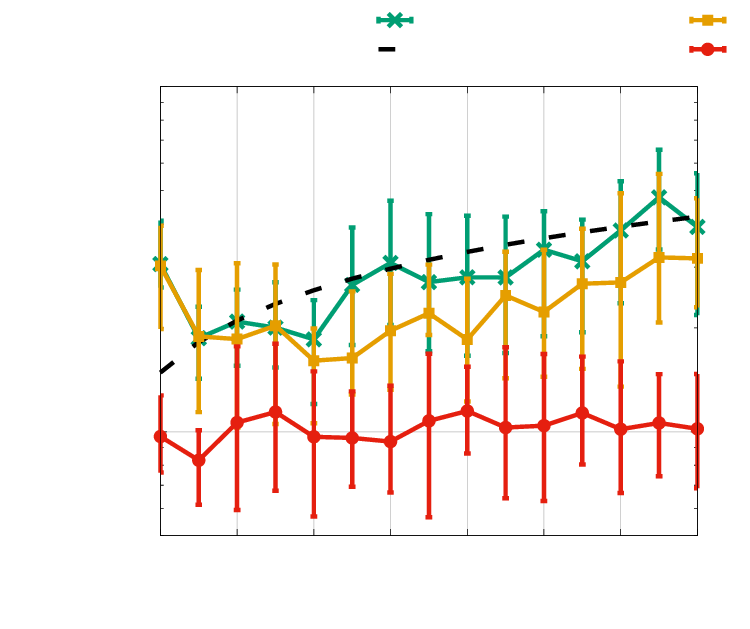}}}\\
 \scriptsize (e) std.~deviation $\sigma = 10^{-4}$, oversampling $\np = 4\nd$ & \scriptsize \hspace*{-0.5cm}(f) std.~deviation $\sigma = 10^{-4}$, oversampling $\np = 8\nd$\\
 \end{tabular}
 \caption{Diffusion reaction example: Taking more sampling points than basis vectors with \odeime leads to stable reduced models in this example. In contrast, reduced models based on \qdeim and \odeimd exhibit instabilities in the sense of Section~\ref{sec:Prelim:Instability}, which is indicated by the growth of the error with the rate $\sqrt{\nd}$. Error bars show the minimum and maximum error over $10$ replicates.}
 \label{fig:NumRes:Laplace:Sigma}
 \end{figure}

 We now run the reduced models for the test parameters and perturb the nonlinear function evaluations $f$ with zero-mean Gaussian noise and standard deviation $\sigma > 0$. We repeat this process $k = 10$ times and collect the states of a reduced model corresponding to the test parameters as columns in $\tilde{\bfX}_i^{(\textsf{S})} \in \mathbb{R}^{\nr \times 81}$ for $i = 1, \dots, k$. Then, the averaged relative state error
 \begin{equation}
 \sum_{i = 1}^k
 \frac{\|\bfX - \bfV \tilde{\bfX}_i^{(\textsf{S})}\|_F}{\|\bfX\|_F}
 \label{eq:NumRes:Laplace:Error}
 \end{equation}
 is reported in the following for each reduced model. Figure~\ref{fig:NLExpAll} compares the error \eqref{eq:NumRes:Laplace:Error} for reduced models based on \qdeim, \odeimd, \odeiml, \odeimr, and \odeime. The standard deviation of the noise is $\sigma = 10^{-2}$ and the oversampling factor is 4, i.e., $\np = 4\nd$. The growth of the error \eqref{eq:NumRes:Laplace:Error} with rate $\sqrt{\nd}$ can be observed for \qdeim in Figure~\ref{fig:NLExpAll}. Similarly, the reduced models based on \odeimd seem unstable because the corresponding errors grow with a rate of $\sqrt{\nd}$ too. In contrast, \odeime and \odeimr give stable reduced models, where the error does not increase with the dimension $\nd$ of the reduced space spanned by the columns of $\bfU$. The curves plotted in Figure~\ref{fig:NLExpAll} are shown in Figure~\ref{fig:NLExpAllWithBars} with error bars that indicate the minimum and maximum error over the $k = 10$ replicates. \cool{The sampling points selected with \odeiml lead to models with poor performance in this example even though the growth of the error with rate $\sqrt{\nd}$ cannot be observed in the plot in Figure~\ref{fig:NLExpAll}. However, the error bars shown in Figure~\ref{fig:NLExpAllWithBars} for \odeiml are larger than for the other sampling algorithms, which indicates that there is strong variability in the approximation error achieved with \odeiml in this example. The strong variability with respect to accuracy of the selected points might hide the growth of the error.} Figure~\ref{fig:NumRes:Laplace:Sigma} compares \odeime with \qdeim and \odeimd for oversampling factors $\np/\nd = 4$ and $\np/\nd = 8$ and standard deviations $\sigma \in \{10^{-2}, 10^{-3}, 10^{-4}\}$. The error bars show the minimum and maximum error over the $k = 10$ replicates. In all cases, \odeime leads to a stable reduced model in the sense that the error does not grow with the dimension $\nd$ of the \deim space, whereas \qdeim and \odeimd show unstable behavior and a growth of the error with rate $\sqrt{\nd}$.

 \begin{figure}
 \centering
 \begin{tabular}{cc}
 {\LARGE\resizebox{0.46\columnwidth}{!}{\input{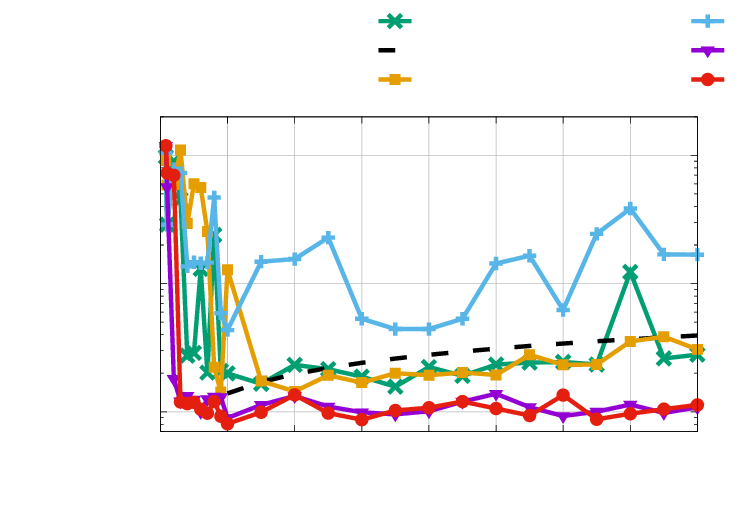}}} & {\LARGE\resizebox{0.46\columnwidth}{!}{\input{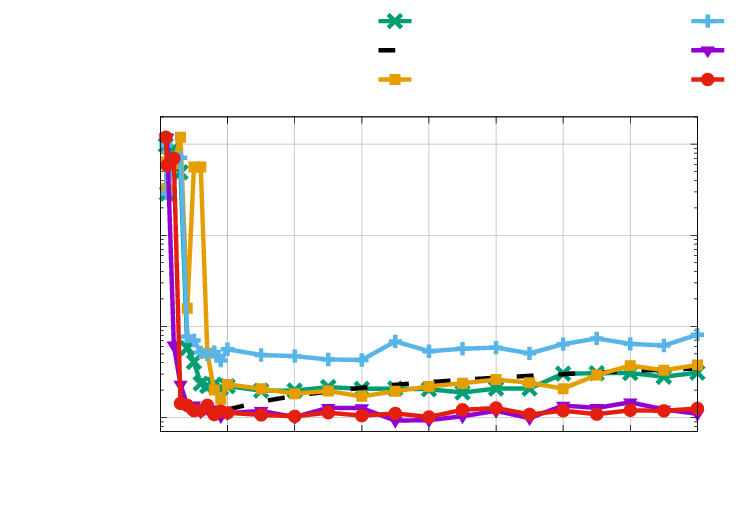}}}\\
 \scriptsize (a) std.~deviation $\sigma = 10^{-1}$ & \scriptsize (b) std.~deviation $\sigma = 10^{-2}$\\
 {\LARGE\resizebox{0.46\columnwidth}{!}{\input{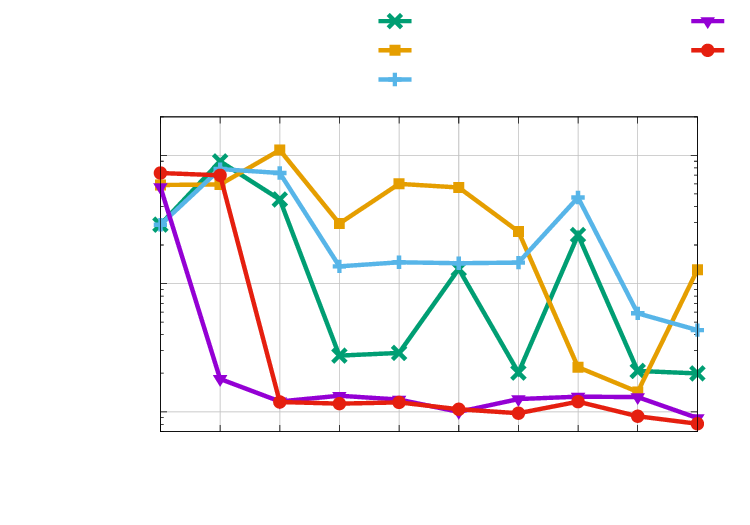}}} & {\LARGE\resizebox{0.46\columnwidth}{!}{\input{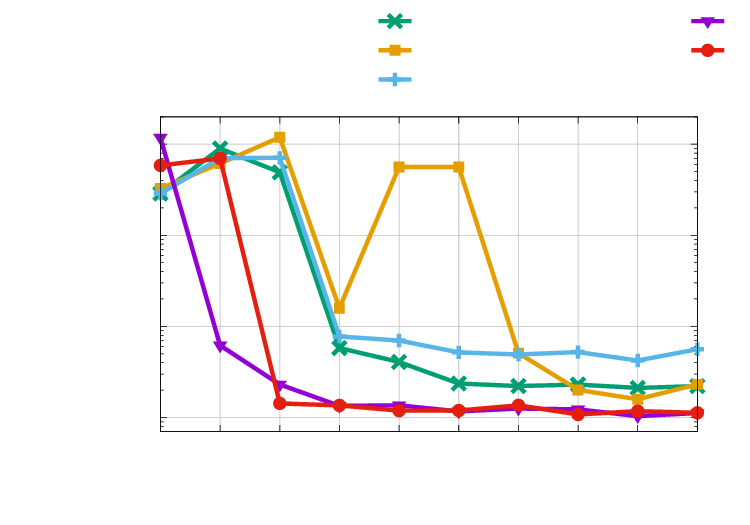}}}\\
 \scriptsize (c) detail of plot (a) & \scriptsize (d) detail of plot (b)\\
 \end{tabular}
 \caption{\cool{Diffusion reaction example: The plots show that \odeime and \odeimr achieve smaller errors than \qdeim also for a lower POD dimension $\nd = 9$ in this example. Notice in the detail of plot (a) shown in (c) that \odeimr and \odeime show less oscillatory error behavior than \qdeim, \odeiml, and \odeimd for DEIM dimensions $\nr$ between 5 and 50.}}
 \label{fig:NLExpLowPODDim}
 \end{figure}

 \cool{Consider now Figure~\ref{fig:NLExpLowPODDim} that shows results for POD dimension $\nr = 9$, which is lower than dimension $\nr = 50$ used previously. The POD space of dimension $\nr = 9$ preserves about 99.9\% of the energy, a typical threshold used in model reduction; cf.~\cite[Section~3.1.1]{BGWSIREV}. Note that the energy is $\sum_{i = 1}^{\nr} \zeta_i^2 \big/ \sum_{i = 1}^{\nh} \zeta_i^2$ where $\zeta_1, \dots, \zeta_{\nh}$ are the singular values of the snapshot matrix in descending order. The standard deviation of the noise is set to $\sigma = 10^{-1}$ and $\sigma = 10^{-2}$, respectively. The mean over 10 runs is shown in Figure~\ref{fig:NLExpLowPODDim}. Similar behavior in terms of error as for higher POD dimensions is observed. Plot (c) in Figure~\ref{fig:NLExpLowPODDim} shows a detail of (a) and indicates that the approximations based on \odeimr and \odeime have less oscillatory error than approximations obtained with \qdeim, \odeiml, and \odeimd for DEIM dimensions $\nd$ between 5 and 50 in this example.}

 \begin{remark}
 \cool{
 We comment on the problem setup: In this example, we considered nonlinear function evaluations that are perturbed with noise. We might encounter such a situation if, for example, parameters of the nonlinear function first need to be estimated from data via a Bayesian approach that introduces noise into the function evaluations used in the reduced model. Our analysis does not cover deterministic approximation errors stemming from, e.g., relaxed tolerances of iterative solvers, and thus it remains future work to show if our analysis applies to such general types of noise as well.}
 \end{remark}

 \section{Conclusions}
 Empirical interpolation is widely used for approximating nonlinear terms in reduced models and for recovering state fields from few spatial measurements; however, stability issues have been observed in presence of noise and other perturbations. Our probabilistic analysis shows that the particular instability that arises due to perturbations such as noise can be provably avoided by employing \gpod and taking more sampling points than dimensions of the reduced space. Numerical results demonstrated that instabilities in \deim can lead to a loss of accuracy in the reduced model outputs and that randomized and deterministic sampling strategies together with \gpod give stabler approximations.

\section*{Acknowledgements}
The authors thank Karthik Duraisamy (University of Michigan), Cheng Huang (University of Mi\-chi\-gan), and David Xu (University of Michigan) for providing the snapshots corresponding to the single-injector combustion process discussed in \Cref{sec:NumRes:Combustion}.

\bibliography{bibDEIM}
\bibliographystyle{siam}

\appendix

\section{Additional listings and figures}

\begin{lstlisting}[style=Matlab-editor,caption={Selecting interpolation points with QDEIM \cite{drmac-gugercin-DEIM-2016} (Matlab code)},label={alg:qdeim}]
function [ q ] = qdeim( U, m )
n = size(U, 2);
[~, ~, q] = qr(U', 'vector');
q = q(1:n)';
end
\end{lstlisting}

\begin{figure}[h]
\centering
{\LARGE\resizebox{0.90\columnwidth}{!}{\input{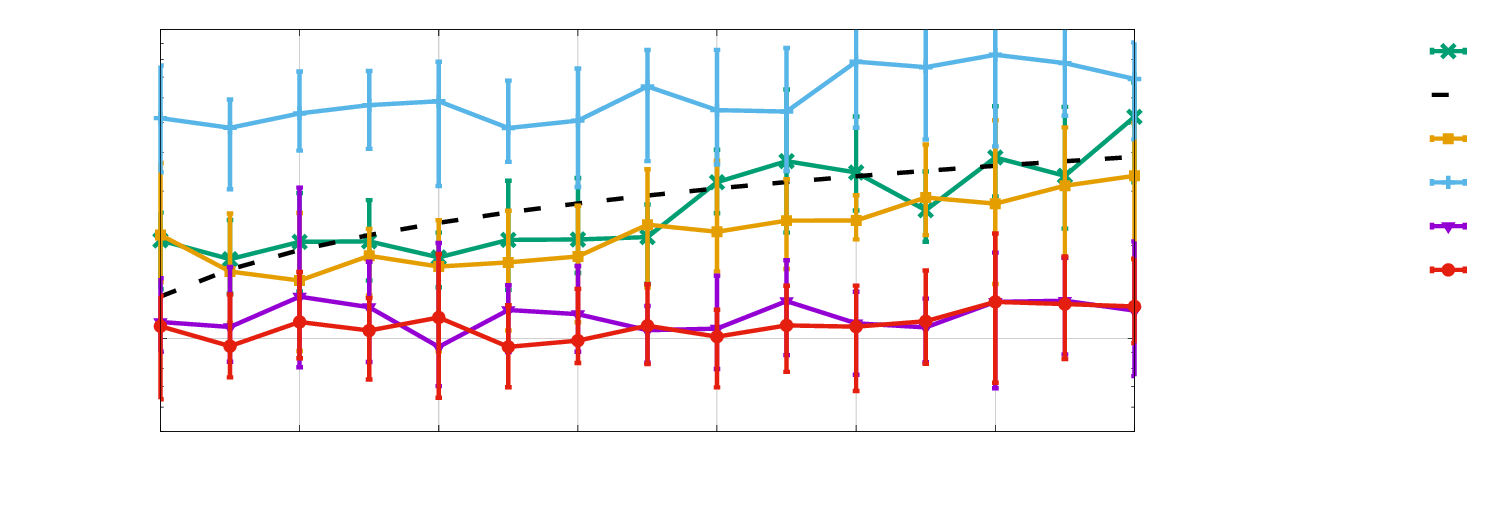}}}\\
\caption{Diffusion reaction example: Approximating the nonlinear terms in this example with \qdeim leads to unstable behavior, which is indicated in this plot with a growth with rate $\sqrt{\nd}$. \gpod with more sampling points than basis vectors leads to stable reduced models with sampling strategies \odeime and \odeimr in this example.  Standard deviation of noise is $\sigma = 10^{-2}$ and oversampling factor is $\np/\nd = 4$.}
\label{fig:NLExpAllWithBars}
\end{figure}

\end{document}